# Polynomial in Non-Commutative Algebra

Aleks Kleyn

ABSTRACT. I considered definition and properties of polynomial in non-commutative algebra. There exists polynomial which has finite, infinite or empty set of roots. For instance, the polynomial
$$p_1(x) = ix - xi - 1$$
have no root and the polynomial
$$p_k(x) = ix - xi - k$$
has the set of roots
$$x = C_1 + C_2 i + \frac{1}{2} j$$
I considered division of polynomials with remainder.

## Contents

## 1. Linear Map

**Convention 1.1.** *If the map*
$$f : A \to A$$
*of $D$-algebra $A$ is linear map, then I use notation*
$$(1.1) \qquad f \circ a = f(a)$$
*for image of the map $f$.* □

Aleks_Kleyn@MailAPS.org.
  http://AleksKleyn.dyndns-home.com:4080/    http://arxiv.org/a/kleyn_a_1.
  http://AleksKleyn.blogspot.com/.

**Definition 1.2.** *Linear map*
$$f : A \to A$$
*of $D$-algebra $A$ satisfies to equalities*

(1.2) $$f \circ (a + b) = f \circ a + f \circ b$$

(1.3) $$f \circ (da) = d(f \circ a)$$

$$a, b \in A \quad d \in D$$

*Let us denote* $\mathcal{L}(D; A \to A)$ *set of linear maps of $D$-algebra $A$.* □

Let $H$ be quaternion algebra.

**Theorem 1.3.** *Let product in algebra $H \otimes H$ be defined according to rule*

(1.4) $$(p_0 \otimes p_1) \circ (q_0 \otimes q_1) = (p_0 q_0) \otimes (q_1 p_1)$$

*A representation*

(1.5) $$h : H \otimes H \dashrightarrow \mathcal{L}(R; H \to H) \quad h(p) : g \to p \circ g$$

*of $R$-algebra $H \otimes H$ in Abelian group $\mathcal{L}(R; H \to H)$ defined by the equality*

(1.6) $$(a \otimes b) \circ g = agb \quad a, b \in H \quad g \in \mathcal{L}(R; H \to H)$$

*is left $H \otimes H$-module. If we put $g = E$ where $E \in \mathcal{L}(R; H \to H)$ is identity map, then we indentify the linear map $f : H \to H$ of quaternion algebra and tensor*

(1.7) $$f = f_{s \cdot 0} \otimes f_{s \cdot 1} \in H \otimes H$$

*using the following equality*

(1.8) $$f \circ a = (f_{s \cdot 0} \otimes f_{s \cdot 1}) \circ a = f_{s \cdot 0} \, a \, f_{s \cdot 1}$$

PROOF. The theorem follows from the statement that the map
$$E \circ x = x$$
generates any linear map of quaternion algebra. □

**Definition 1.4.** *Expression $f_{s \cdot p}$, $p = 0, 1$, in equality (1.7) is called **component of linear map** $f$.* □

$$\boxed{(1.7) \quad f = f_{s \cdot 0} \otimes f_{s \cdot 1} \in H \otimes H}$$

**Definition 1.5.** *The polylinear map of $D$-algebra $A$*
$$f : A^n \to A$$
*satisfies to equalities*

(1.9) $$f \circ (a_1, ..., a_i + b_i, ..., a_n) = f \circ (a_1, ..., a_i, ..., a_n) + f \circ (a_1, ..., b_i, ..., a_n)$$

(1.10) $$f \circ (a_1, ..., p a_i, ..., a_n) = p f \circ (a_1, ..., a_i, ..., a_n)$$

$$1 \le i \le n \quad a_i, b_i \in A \quad p \in D$$

*Let us denote* $\mathcal{L}(D; A^n \to A)$ *set of $n$-linear maps of $D$-algebra $A$.* □

**Theorem 1.6.** *A representation*

$$h : H^{n+1\otimes} \times S_n \dashrightarrow \mathcal{LA}(R; H^n \to H)$$

*of algebra $H^{n+1\otimes}$ in module $\mathcal{LA}(R; H^n \to H)$ defined by the equality*

(1.11)
$$(a_0 \otimes ... \otimes a_n, \sigma) \circ (f_1 \otimes ... \otimes f_n) = a_0 \sigma(f_1) a_1 ... a_{n-1} \sigma(f_n) a_n$$
$$a_0, ..., a_n \in H \quad \sigma \in S_n \quad f_1, ..., f_n \in \mathcal{L}(R; H \to H)$$

*allows us to identify tensor $d \in A^{n+1\otimes}$ and transposition $\sigma \in S^n$ with map*

(1.12) $\quad (d, \sigma) \circ (f_1, ..., f_n) \quad f_i = \delta \in \mathcal{L}(R; H \to H)$

*where $\delta \in \mathcal{L}(R; H \to H)$ is identity map.*

## 2. Polynomial

**Definition 2.1.** *Let*
$$f : A^n \to A$$
*be polylinear map of D-module A. The map*
$$p : A \to A$$
*defined by the equality*

(2.1) $\quad p \circ x^n = f \circ (x_1, ..., x_n) \quad x_1 = ... = x_n = x$

*is called* **homogeneous polynomial** *of power $n$.* □

**Theorem 2.2.** *The set $A_n[x]$ of homogeneous polynomials of power $n$ is Abelian group.*

Proof. The theorem follows from the statement that sum of $n$-linear maps is $n$-linear map. □

**Theorem 2.3.** *Let $A$ be D-algebra. Then for any homogeneous polynomial*
$$p : A \to A$$
*of power $n$ there exists tensor $p' \in A^{n+1\otimes}$ such that*

(2.2) $\quad p \circ x^n = p' \circ x^n = p' \circ (x_1, ..., x_n) \quad x_1 = ... = x_n = x$

Proof. The theorem follows the theorem 1.6. Hereinafter we will identify homogeneous polynomial $p$ and corresponding tensor $p'$. □

**Definition 2.4.** *If the tensor $p$ has form $p = p_0 \otimes ... \otimes p_n$, then polynomial $p \circ x^n$ is called monomial.* □

**Theorem 2.5.** *Let $p_k(x)$ be* **monomial of power** $k$ *over associative D-algebra $A$. Then*

2.5.1: *Monomial of power* 0 *has form* $p_0(x) = a_0$, $a_0 \in A$.
2.5.2: *If* $k > 0$, *then*

$$p_k(x) = p_{k-1}(x)xa_k \tag{2.3}$$

*where* $a_k \in A$.

PROOF. We prove the theorem by induction over power $n$ of monomial.

Let $n = 0$. We get the statement 2.5.1 since monomial $p_0(x)$ is constant.

Let $n = k$. Last factor of monomial $p_k(x)$ is either $a_k \in A$, or has form $x^l$, $l \geq 1$. In the later case we assume $a_k = 1$. Factor preceding $a_k$ has form $x^l$, $l \geq 1$. We can represent this factor as $x^{l-1}x$. Therefore, we proved the statement. □

**Remark 2.6.** *In the theorem 2.5, I considered recursive representation of a monomial in associative D-algebra. Since product is independent of the way in which brackets are placed, recursive representation of a monomial is not unique in associative D-algebra. For instance, I can use any of the following forms*

$$\begin{aligned}ax^2bxcx^3d &= ax(xbxcx^3) = (ax)x(bxcx^3d) = (ax^2b)x(cx^3d)\\&= (ax^2bxc)x(x^2d) = (ax^2bxcx)x(xd) = (ax^2bxcx^2)xd\end{aligned} \tag{2.4}$$

*to represent monomial* $ax^2bxcx^3d$. *I chose the equality* (2.3) *as the most simple for algorithm of division of polynomials.* □

**Definition 2.7.** *We denote*

$$A[x] = \bigoplus_{n=0}^{\infty} A_n[x]$$

*direct sum*[2.1] *of Abelian groups* $A_n[x]$. *An element* $p(x)$ *of Abelian group* $A[x]$ *is called* **polynomial** *over D-algebra* $A$. □

Therefore, we can present polynomial of power $n$ in the following form

$$p(x) = p_0 + p_1 \circ x + ... + p_n \circ x^n \quad p_i \in A^{i+1\otimes} \quad i = 0, ..., n \tag{2.5}$$

## 3. OPERATIONS WITH POLYNOMIALS

**Definition 3.1.** *Let*

$$p(x) = p_0 + p_1 \circ x + ... + p_n \circ x^n \quad p_i \in A^{i+1\otimes} \quad i = 0, ..., n$$

$$r(x) = r_0 + r_1 \circ x + ... + r_n \circ x^n \quad r_i \in A^{i+1\otimes} \quad i = 0, ..., n$$

*be polynomials.*[3.1] *We introduce the sum of polynomials* $p$ *and* $r$ *by the equality*

$$(p + r)(x) = p_0 + r_0 + (p_1 + r_1) \circ x + ... + (p_n + r_n) \circ x^n \tag{3.1}$$

□

---

[2.1] See the definition of direct sum of Abelian groups in [1], pages 36, 37. On the same page, Lang proves the existence of direct sum of Abelian groups.

[3.1] If the coefficient, let's say $p_i$, is absent, then we assume $p_i = 0$.

**Definition 3.2.** *Bilinear map*

(3.2) $$*: A^{n\otimes} \times A^{m\otimes} \to A^{n+m-1\otimes}$$

*is defined by the equality*

(3.3) $$(a_1 \otimes ... \otimes a_n) * (b_1 \otimes ... \otimes b_n) = a_1 \otimes ... \otimes a_{n-1} \otimes a_n b_1 \otimes b_2 \otimes ... \otimes b_n$$

□

**Theorem 3.3.** *For any tensors* $a \in A^{n+1\otimes}$, $b \in A^{m+1\otimes}$, *product of homogeneous polynomials* $a \circ x^n$, $b \circ x^m$ *is defined by the equality*

(3.4) $$(a * b) \circ x^{n+m} = (a \circ x^n)(b \circ x^m)$$

## 4. Linear Equation

Let $\bar{\bar{e}}$ be the basis of finite dumentional algebra $A$ over field $F$ and $C_{ij}^k$ be structural constants of algebra $A$ relative to the basis $\bar{\bar{e}}$.

Consider the linear equation

(4.1) $$a \circ x = b$$

where $a = a_{s\cdot 0} \otimes a_{s\cdot 1} \in A^{2\otimes}$. According to the theorem [3]-6.4.1, we can write the equation (4.1) in standard form

(4.2) $$a^{ij} e_i x e_j = b$$

(4.3) $$a^{ij} = a_{s\cdot 0}{}^i a_{s\cdot 1}{}^j \quad a_{s\cdot 0} = a_{s\cdot 0}{}^i e_i \quad a_{s\cdot 1} = a_{s\cdot 1}{}^i e_i$$

According to the theorem [3]-6.4.5, equation (4.2) is equivalent to equation

(4.4) $$a_i^j x^i = b^j$$

(4.5) $$a_i^j = a^{kr} C_{ki}^p C_{pr}^j \quad x = x^i e_i \quad b = b^i e_i$$

According to the theory of linear equations over field, if determinant

(4.6) $$\det \|a_i^j\| \neq 0$$

then equation (4.1) has only one solution.

**Definition 4.1.** *The tensor* $a \in A^{2\otimes}$ *is called* **nonsingular tensor** *if this tensor satisfies to condition* (4.6). □

**Theorem 4.2.** *Let* $a \in A^{2\otimes}$ *be nonsingular tensor. If we consider the equation* (4.2) *as transformation of algebra $A$, then we can write the inverse transformation in form*

(4.7) $$x = c^{pq} e_p b e_q$$

*where components $c^{pq}$ satisfy to equation*

$$\delta_0^r \delta_0^s = a^{ij} c^{pq} C_{ip}^r C_{qj}^s \tag{4.8}$$

PROOF. The theorem follows from the theorem [2]-6.7. □

| (4.1) $a \circ x = b$ | (4.2) $a^{ij} e_i x e_j = b$ |

[2] Aleks Kleyn, Polynomial over Associative $D$-Algebra, eprint arXiv:1302.7204 (2013)

**Definition 4.3.** *Let $a \in A^{2\otimes}$ be nonsingular tensor. The tensor*

$$a^{-1} = c^{pq} e_p \otimes e_q \tag{4.9}$$

*is called* **tensor inverse to tensor** $a$. □

**Theorem 4.4.** *Let $a \in A^{2\otimes}$ be nonsingular tensor. Then linear equation*

$$a \circ x = b \tag{4.10}$$

*has unique solution*

$$x = a^{-1} \circ b \tag{4.11}$$

*If $a$ is singular tensor, then the equation (4.10) has solution only when $b \in \ker a$. In this case the equation (4.10) has infinitely many solutions.*

| (4.10) $a \circ x = b$ |

**Example 4.5.** *The tensor*

$$i \otimes 1 - 1 \otimes i \in H^{2\otimes}$$

*corresponds to the matrix*

$$\begin{pmatrix} 0 & -1 & 0 & 0 \\ 1 & 0 & 0 & 0 \\ 0 & 0 & 0 & -1 \\ 0 & 0 & 1 & 0 \end{pmatrix} - \begin{pmatrix} 0 & -1 & 0 & 0 \\ 1 & 0 & 0 & 0 \\ 0 & 0 & 0 & 1 \\ 0 & 0 & -1 & 0 \end{pmatrix} = \begin{pmatrix} 0 & 0 & 0 & 0 \\ 0 & 0 & 0 & 0 \\ 0 & 0 & 0 & -2 \\ 0 & 0 & 2 & 0 \end{pmatrix} \tag{4.12}$$

*Therefore, the polynomial*

$$p_1(x) = ix - xi - 1 \tag{4.13}$$

*have no root and the polynomial*

$$p_k(x) = ix - xi - k \tag{4.14}$$

*has the set of roots*

$$x = C_1 + C_2 i + \frac{1}{2} j \tag{4.15}$$

□

## 5. Left-sided Polynomial

[6] Paul M. Cohn, Skew Fields, Cambridge University Press, 1995

In the book [6] on the page 48, professor Cohn wrote that the study of polynomials over skew field is difficult problem. For this reason professor Ore considered skew right polynomial ring [5.1] where Ore introduced the transformation of monomial

(5.1) $$at = ta^\alpha + a^\delta$$

The map

$$\alpha : a \in A \to a^\alpha \in A$$

is endomorphism of $D$-algebra $A$ and the map

$$\delta : a \in A \to a^\delta \in A$$

is $\alpha$-derivation of $D$-algebra $A$

(5.2) $$(a+b)^\delta = a^\delta + b^\delta \quad (ab)^\delta = a^\delta b^\alpha + ab^\delta$$

In such case all coefficients of polynomial can be written on the right. Such polynomials are called right-sided polynomials.

The same way we may consider the transformation of monomial

(5.3) $$ta = a^\alpha t + a^\delta$$

In such case all coefficients of polynomial can be written on the left. Such polynomials are called left-sided polynomials. [5.2]

The point of view proposed by Ore allows us to get qualitative picture of the theory of polynomials in non-commutative algebra. However, there exist algebras where this point of view is not complete.

**Example 5.1.** *If there exist transformation (5.1) in $D$-algebra $A$, then we can map any polynomial $p$ of power $1$ into right-sided polynomial $p'$ of power $1$ such way that*

$$p(t) = p'(t)$$

*The polynomial $p'$ may have the following form*

(5.4) $$p'(t) = ta_1 + a_0 \quad a_1 \neq 0$$

(5.5) $$p'(t) = t\,0 + a_0 \quad a_0 \neq 0$$

(5.6) $$p'(t) = t\,0 + 0$$

*The polynomial (5.4) has unique solution. The polynomial (5.5) has no solution. Any $A$-number is root of the polynomial (5.6).*

*It is evident that we missed the case when the polynomial $p$ has infinitely many roots, but not any $A$-number is the root of the polynomial $p$. Therefore the question about the existence of transformation (5.1) in quaternion algebra arises.* □

---

[5.1] See also the definition on pages [6]-(48 - 49).
[5.2] See also the definition on the page [5]-3.

8                                       Aleks Kleyn

**Example 5.2.** *Let*

$$p(x) = \sum_{i=0}^{n} p_i x^i$$

*and*

$$r(x) = \sum_{j=0}^{m} r_j x^j$$

*be left-sided polynomials. We introduce product of polynomials p and r by the equality*[5.3]

(5.7) $$p(x) * r(x) = (p * r)(x) = \sum_{j=0}^{m+n} (\sum_{i=0}^{j} p_i r_{j-i}) x^j$$

*Consider polynomials*

(5.8) $$L(x) = x - i$$

(5.9) $$R(x) = x - j$$

*The equality*

(5.10) $$P(x) = (L * R)(x) = x^2 - (i+j)x + k$$

*follows from equalities* (5.7), (5.8), (5.9). *Let*

$$h = R(i) = i - j$$

$$\tilde{x} = hxh^{-1} = (i-j)x(i-j)^{-1} = \frac{1}{2}(i-j)x(j-i)$$

(5.11) $$L(\tilde{x}) = \frac{1}{2}(i-j)x(j-i) - i$$

*According to the theorem [5]-1, item (ii), the polynomial*

(5.12) $$P_1(x) = L(\tilde{x})R(x) = (\frac{1}{2}(i-j)x(j-i) - i)(x - j)$$

*coincides with the polynomial* $P(x)$, $P_1(x) = P(x)$

(5.13) $$x^2 - (i+j)x + k = (\frac{1}{2}(i-j)x(j-i) - i)(x - j)$$

*However, the equality* (5.13) *is not true. Let* $x = i + j$. *Equalities*

(5.14) $$(i+j)^2 - (i+j)(i+j) + k = (\frac{1}{2}(i-j)(i+j)(j-i) - i)(i+j-j)$$

(5.15) $$\begin{aligned} k &= (\frac{1}{2}(i^2 + ij - ji - j^2)(j-i) - i)i \\ &= (k(j-i) - i)i = (-i - j - i)i \end{aligned}$$

*follow from the equality* (5.13).                                                    □

6. Square Root

---

[5.3] See also the definition on the page [5]-3.

**Definition 6.1.** *The root* $x = \sqrt{a}$ *of the equation*

(6.1) $$x^2 = a$$

*in D-algebra A is called* **square root** *of A-number a.* □

**Theorem 6.2.** *Let H be quaternion algebra and a be H-number.*

6.2.1: *Since* $\operatorname{Re}\sqrt{a} \neq 0$, *then the equation*

$$x^2 = a$$

*has roots* $x = x_1$, $x = x_2$ *such that*

(6.2) $$x_2 = -x_1$$

6.2.2: *Since* $a = 0$, *then the equation*

$$x^2 = a$$

*has root* $x = 0$ *with multiplicity* $2$.

6.2.3: *Since conditions* 6.2.1, 6.2.2 *are not true, then the equation*

$$x^2 = a$$

*has infinitly many roots such that*

(6.3) $$x \in \operatorname{Im} H \quad a \in \operatorname{Re} H \quad |x| = \sqrt{-a}$$

PROOF. The theorem follows from the theorem [4]-8.5. □

**Theorem 6.3.** *Let quaternion a be a a root of the polynomial*

(6.4) $$r(x) = x^2 + 1$$

*Then*

(6.5) $$\sqrt{a} = \pm \frac{1}{\sqrt{2}}(1 + a)$$

PROOF. The equality

(6.6) $$(a+1)^2 = a^2 + 2a + 1 = 2a$$

follows from the statement that $a$ is a root of the polynomial (6.4). The equality (6.5) follows from the equality (6.6) and from the statement 6.2.1. □

## 7. Polynomial with Given Roots

In commutative algebra, if I have two roots of polynomial of second order, then this polynomial is unique. The case is different in non-commutative algebra.

**Theorem 7.1.** *Let* $x_1$, $x_2 \in A$,

(7.1) $$x_1 x_2 \neq x_2 x_1$$

*Then*

(7.2) $$p_{12}(x) = (x - x_1)(x - x_2) = x^2 - x_1 x - x x_2 + x_1 x_2$$

$$(7.3) \qquad p_{21}(x) = (x - x_2)(x - x_1) = x^2 - x_2 x - x x_1 + x_2 x_1$$

are unequal polynomials

$$(7.4) \qquad p_{12} \neq p_{21}$$

PROOF. Let

$$(7.5) \qquad x = x_1 + x_2$$

The equality

$$(7.6) \qquad p_{12}(x_1 + x_2) = (x_1 + x_2 - x_1)(x_1 + x_2 - x_2) = x_2 x_1$$

follows from equalities (7.2), (7.5). The equality

$$(7.7) \qquad p_{21}(x_1 + x_2) = (x_1 + x_2 - x_2)(x_1 + x_2 - x_1) = x_1 x_2$$

follows from equalities (7.3), (7.5). The statement (7.4) follows from equalities (7.6), (7.7) and from the statement (7.1). □

**Theorem 7.2.** *Let $a_1$, $a_2 \in A \otimes A$. Then the polynomial*

$$(7.8) \qquad p(x) = a_1 \circ p_{12}(x) + a_2 \circ p_{21}(x)$$

*has roots $x = x_1$, $x = x_2$.*

PROOF. The theorem follows from the statement that polynomials $p_{12}$, $p_{21}$ have roots $x = x_1$, $x = x_2$. □

**Question 7.3.** *Does polynomial (7.8) have roots which are different from values $x = x_1$, $x = x_2$?* □

**Question 7.4.** *Are there other polynomials that have roots $x = x_1$, $x = x_2$?* □

It looks like the answer to the question 7.4 is positive.

**Example 7.5.** *Consider the set of polynomials*

$$(7.9) \qquad p(x) = a_1 \circ ((x - i)(x - j)) + a_2 \circ ((x - j)(x - i))$$

*which, according to the theorem 7.2, have roots $x = i$, $x = j$. The polynomial*

$$(7.10) \qquad r(x) = x^2 + 1$$

*also has roots $x = i$, $x = j$.* □

**Question 7.6.** *Does polynomial (7.10) belong to the set of polynomials (7.9)?* □

The theorem 7.7 answers the question 7.6.

**Theorem 7.7.** *There is no tensors $a_1$, $a_2 \in A \otimes A$ such that*

$$(7.11) \qquad a_1 \circ (x^2 - ix - xj + k) + a_2 \circ (x^2 - jx - xi - k) = x^2 + 1$$

PROOF. Let the theorem be not true.

**Statement 7.8.** *Let there exist tensors $a_1$, $a_2 \in A \otimes A$ such that the equalitiy (7.11) is true.* ⊙

Equalities

$$a_1 \circ x^2 + a_2 \circ x^2 = x^2 \tag{7.12}$$

$$a_1 \circ (ix) + a_1 \circ (xj) + a_2 \circ (jx) + a_2 \circ (xi) = 0 \tag{7.13}$$

$$a_1 \circ k - a_2 \circ k = 1 \tag{7.14}$$

follow from the equality (7.11).

Let $x = 1$. The equality

$$a_1 \circ i + a_1 \circ j + a_2 \circ j + a_2 \circ i = 0 \tag{7.15}$$

follow from equalities (7.13).

According to the theorem 6.3, there exists x such that $x^2 = i$. The equality

$$a_1 \circ i + a_2 \circ i = i \tag{7.16}$$

follows from the equality (7.12).

According to the theorem 6.3, there exists x such that $x^2 = j$. The equality

$$a_1 \circ j + a_2 \circ j = j \tag{7.17}$$

follows from the equality (7.12).

The equality

$$a_1 \circ i + a_2 \circ i + a_1 \circ j + a_2 \circ j = i + j \tag{7.18}$$

follows from equalities (7.16), (7.17) and contradicts the equality (7.15). Therefore, the statement 7.8 is not true and we proved the theorem. □

## 8. Division with Remainder

**Definition 8.1.** *A-number $a$ is left divisor of A-number $b$, if there exists A-number $c$ such that*

$$ac = b \tag{8.1}$$

□

**Definition 8.2.** *A-number $a$ is right divisor of A-number $b$, if there exists A-number $c$ such that*

$$ca = b \tag{8.2}$$

□

It is evident that there is symmetry between definitions 8.1 and 8.2. The difference between left and right divisors is also evident since the product is noncommutative. However we can consider a definition generalizing definitions 8.1 and 8.2.

**Definition 8.3.** *A-number $a$ is divisor of A-number $b$, if there exists $A \otimes A$-number $c$ such that*

$$c \circ a = b \tag{8.3}$$

*$A \otimes A$-number $c$ is called* **quotient** *of A-number $b$ divided by A-number $a$.* □

**Definition 8.4.** *Let division in the D-algebra $A$ is not always defined. A-**number $a$ divides A-number $b$ with remainder**, if the following equation is true*

$$c \circ a + f = b \tag{8.4}$$

*$A \otimes A$-number $c$ is called* **quotient** *of A-number $b$ divided by A-number $a$. A-number $f$ is called* **remainder of the division** *of A-number $b$ by A-number $a$.* □

**Theorem 8.5.** *Let*

$$p(x) = x - a$$

*be polynomial of power $1$. Let*

$$r(x) = r_0 + r_1 \circ x + ... + r_k \circ x^k$$

*be polynomial of power $k > 0$. Then*

$$r(x) = s_0 + (q_0 + q_1(x) + ... + q_{k-1}(x)) \circ p(x) \tag{8.5}$$

*where $q_i \in A_i[x] \otimes A$, $i = 0, ..., k-1$, is homogeneous polynomial of power $i$.*

Proof. Let

$$r_k = r_{k.0.s} \otimes ... \otimes r_{k.k.s} \tag{8.6}$$

According to the theorem 2.5,

$$\begin{aligned}r_k \circ x^k &= ((r_{k.0.s} \otimes ... \otimes r_{k.k-1.s}) \circ x^{k-1})xr_{k.k.s}\\ &= (((r_{k.0.s} \otimes ... \otimes r_{k.k-1.s}) \circ x^{k-1}) \otimes r_{k.k.s}) \circ x\end{aligned} \tag{8.7}$$

Let

$$q_{k-1}(x) = ((r_{k.0.s} \otimes ... \otimes r_{k.k-1.s}) \circ x^{k-1}) \otimes r_{k.k.s} \tag{8.8}$$

To prove the theorem, it is enough to note that power of the polynomial

$$s_{k-1}(x) = r(x) - q_{k-1}(x) \circ (x - a) \tag{8.9}$$

is less than $k$. The equality

$$r(x) = s_{k-1}(x) + q_{k-1}(x) \circ (x - a) \tag{8.10}$$

follows from the equality (8.9). We will apply algorithm considered in proof to the polynomial $s_{k-1}(x)$. After a finite number of steps, we will get the polynomial $s_0$ of power $0$. The equality (8.5) follows from the equality (8.10). □

The algorithm considered in the proof of the theorem is called standard algorithm of division of polynomials.

## 9. Examples of Division of Polynomials

In following examples, we use standard algorithm of division of polynomials.

**Example 9.1.** *Let*
$$p(x) = x - i$$
(9.1) $$r(x) = (x-j)(x-i) = x^2 - jx - xi - k$$

*According to the theorem 8.5, we set* $q_1(x) = x \otimes 1$. *Then*

(9.2) $$\begin{aligned} s_1(x) &= r(x) - q_1(x) \circ p(x) = r(x) - (x \otimes 1) \circ (x-i) \\ &= x^2 - jx - xi - k - x(x-i) = x^2 - jx - xi - k - x^2 + xi \\ &= -jx - k \end{aligned}$$

*Now we set* $q_0 = -j \otimes 1$. *Then*

(9.3) $$\begin{aligned} s_0 &= s_1(x) - q_0 \circ p(x) = s_1(x) - (-i \otimes 1 - 1 \otimes j + 1 \otimes i) \circ (x-i) \\ &= -ix - xj + k + xi + i(x-i) + (x-i)j - (x-i)i \\ &= -ix - xj + k + xi + ix + 1 + xj - k - xi - 1 = 0 \end{aligned}$$

*The equality*

(9.4) $$\begin{aligned} r(x) &= s_1(x) + (x \otimes 1) \circ (x-i) \\ &= (-j \otimes 1 + x \otimes 1) \circ (x-i) \end{aligned}$$

*follows from equalities* (9.2), (9.3). *We can reduce expression in the equality* (9.4)

(9.5) $$r(x) = ((x-j) \otimes 1) \circ (x-i) = (x-j)(x-i)$$

*We see that the representation of the polynomial* $r(x)$ *in the equality* (9.5) *is the same as the representation of the polynomial* $r(x)$ *in the equality* (9.1). $\square$

The answer (9.5) in the example 9.1 is evident because I intentionally chose divisor $p(x)$ to be equal second factor of the polynomial $r(x)$.

$\boxed{p(x) = x - i}$  $\boxed{r(x) = (x-j)(x-i) = x^2 - jx - xi - k}$

The example 9.2 is more interesting because I considered different order of factors.

**Example 9.2.** *Let*
$$p(x) = x - i$$
(9.6) $$r(x) = (x-i)(x-j) = x^2 - ix - xj + k$$

*According to the theorem 8.5, we set* $q_1(x) = x \otimes 1$. *Then*

(9.7) $$\begin{aligned} s_1(x) &= r(x) - (x \otimes 1) \circ (x-i) \\ &= x^2 - ix - xj + k - x(x-i) = x^2 - ix - xj + k - x^2 + xi \\ &= -ix - xj + k + xi \end{aligned}$$

*Now we set*

(9.8) $$q_0 = -i \otimes 1 - 1 \otimes j + 1 \otimes i$$

*Then*

(9.9) $$\begin{aligned} s_0 &= s_1(x) - q_0 \circ p(x) = s_1(x) - (-j \otimes 1) \circ (x-i) \\ &= -jx - k + j(x-i) = -jx - k + jx + k = 0 \end{aligned}$$

*The equality*

(9.10)
$$r(x) = s_1(x) + (x \otimes 1) \circ (x - i)$$
$$= (-i \otimes 1 - 1 \otimes j + 1 \otimes i + x \otimes 1) \circ (x - i)$$

*follows from equalities* (9.7), (9.9). *We can reduce expression in the equality* (9.10)

(9.11)
$$r(x) = -i(x - i) - (x - i)j + (x - i)i + x(x - i)$$
$$= (x - i)(i - j) + (x - i)(x - i)$$
$$= (x - i)(x - i + i - j)$$
$$= (x - i)(x - j)$$

*We see that the representation of the polynomial $r(x)$ in the equality* (9.11) *is the same as the representation of the polynomial $r(x)$ in the equality* (9.6). □

The quotient of polynomial $r(x)$ divided by polynomial $p(x)$ in the example 9.2 has a more complicated structure than the quotient in the example 9.1. This is because I am trying to write down the divisor $p(x)$ to the right of the quotient, although initially the polynomial $p(x)$ was left factor. Nevertheless, as a result of simple transformations, factorization in the equality (9.11) took the expected form.

$\boxed{p(x) = x - i}$ $\boxed{r(x) = (x - i)(x - j) = x^2 - ix - xj + k}$

In the example 9.3, I consider division of the left-sided polynomial

$$r(x) = x^2 - ix - jx - k$$

over the polynomial

$$p(x) = x - i$$

I know that the polynomial $r(x)$ has solution $x = i$ and want to find another solution.

**Example 9.3.** *Let*

$$p(x) = x - i$$
$$r(x) = x^2 - ix - jx - k$$

*According to the theorem* 8.5, *we set* $q_1(x) = x \otimes 1$. *Then*

(9.12)
$$s_1(x) = r(x) - (x \otimes 1) \circ (x - i)$$
$$= x^2 - ix - jx - k - x(x - i) = x^2 - ix - jx - k - x^2 + xi$$
$$= -ix - jx - k + xi$$

*Now we set*

(9.13)
$$q_0 = -i \otimes 1 - j \otimes 1 + 1 \otimes i$$

*Then*

(9.14)
$$s_0 = s_1(x) - q_0 \circ p(x) = s_1(x) - (-i \otimes 1 - j \otimes 1 + 1 \otimes i) \circ (x - i)$$
$$= -ix - jx - k + xi + i(x - i) + j(x - i) - (x - i)i$$
$$= -ix - jx - k + xi + ix + 1 + jx + k - xi - 1 = 0$$

*The equality*

(9.15)
$$r(x) = s_1(x) + (x \otimes 1) \circ (x - i) = q_0 \circ (x - i) + (x \otimes 1) \circ (x - i)$$
$$= (-i \otimes 1 - j \otimes 1 + 1 \otimes i + x \otimes 1) \circ (x - i)$$

*follows from equalities* (9.12), (9.14). *We can reduce expression in the equality* (9.15)

$$r(x) = -i(x-i) - j(x-i) + (x-i)i + x(x-i) \tag{9.16}$$
$$= (x-i)i + (x-i-j)(x-i)$$

□

**Question 9.4.** *Does the polynomial*

$$r(x) = x^2 - ix - jx - k$$

*have root different from $x = i$?*  □

We verify immediately that $x = i$ is the root of the polynomial $r(x)$. According to the example 9.3, we can represent the polynomial $r(x)$ as

$$r(x) = (1 \otimes i + (x - i - j) \otimes 1) \circ (x - i) \tag{9.17}$$

Therefore, in order for the polynomial $p(x)$ to have the root different from $x = i$, it is necessary that linear map

$$1 \otimes i + (x - i - j) \otimes 1 \tag{9.18}$$

has non-trivial kernel and quaternion $x - i$ belongs to the kernel of the linear map (9.18). The linear map (9.18) has the matrix

$$\begin{pmatrix} 0 & -1 & 0 & 0 \\ 1 & 0 & 0 & 0 \\ 0 & 0 & 0 & 1 \\ 0 & 0 & -1 & 0 \end{pmatrix} + \begin{pmatrix} x^0 & -x^1+1 & -x^2+1 & -x^3 \\ x^1-1 & x^0 & -x^3 & x^2-1 \\ x^2-1 & x^3 & x^0 & -x^1+1 \\ x^3 & -x^2+1 & x^1-1 & x^0 \end{pmatrix}$$
$$\tag{9.19}$$
$$= \begin{pmatrix} x^0 & -x^1 & -x^2+1 & -x^3 \\ x^1 & x^0 & -x^3 & x^2-1 \\ x^2-1 & x^3 & x^0 & -x^1+2 \\ x^3 & -x^2+1 & x^1-2 & x^0 \end{pmatrix}$$

Search of a value of $x$ for which the matrix (9.19) is not singular is not simple task. However I hope to do research in this direction.

In the example 9.5, I consider division of the left-sided polynomial

$$r(x) = x^2 - ix - jx - k \qquad r(j) = -2k$$

over the polynomial

$$p(x) = x - j$$

**Example 9.5.** *Let*

$$p(x) = x - j$$
$$r(x) = x^2 - ix - jx - k$$

*According to the theorem 8.5, we set* $q_1(x) = x \otimes 1$. *Then*

$$s_1(x) = r(x) - (x \otimes 1) \circ (x - i)$$
(9.20)
$$= x^2 - ix - jx - k - x(x - j) = x^2 - ix - jx - k - x^2 + xj$$
$$= -ix - jx - k + xj$$

*Now we set*

(9.21) $$q_0 = -i \otimes 1 - j \otimes 1 + 1 \otimes j$$

*Then*

$$s_0 = s_1(x) - q_0 \circ p(x) = s_1(x) - (-i \otimes 1 - j \otimes 1 + 1 \otimes j) \circ (x - j)$$
(9.22)
$$= -ix - jx - k + xj + i(x - j) + j(x - j) - (x - j)j$$
$$= -ix - jx - k + xj + ix - k + jx + 1 - xj - 1 = -2k$$

*The equality*

$$r(x) = s_1(x) + (x \otimes 1) \circ (x - j)$$
(9.23)
$$= s_0 + q_0 \circ (x - j) + (x \otimes 1) \circ (x - j)$$
$$= -2k + (-i \otimes 1 - j \otimes 1 + 1 \otimes j + x \otimes 1) \circ (x - j)$$

*follows from equalities* (9.20), (9.22). *We can reduce expression in the equality* (9.23)

(9.24) $$r(x) = -2k + (1 \otimes j + (x - i - j) \otimes 1) \circ (x - j)$$

$\square$

## 10. More Questions

In the example 10.1, I considered the polynomial $r(x)$ of power 3 in order to find quotient of the polynomial $r(x)$ divided by two polynomials. First I will find a quotient of the polynomial

$$r(x) = (x - j)(x - k)(x - j - k)$$

divided by the polynomial

$$p(x) = x - k$$

**Example 10.1.** *Let*

$$p(x) = x - k$$

$$r(x) = (x - j)(x - k)(x - j - k) = (x^2 - jx - xk + i)(x - j - k)$$
$$= x^2(x - j - k) - jx(x - j - k) - xk(x - j - k) + i(x - j - k)$$
$$= x^3 - x^2j - x^2k - jx^2 + jxj + jxk - xkx + xkj - x + ix - k + j$$
$$= x^3$$
$$+ (-1 \otimes 1 \otimes j - 1 \otimes 1 \otimes k - j \otimes 1 \otimes 1 - 1 \otimes k \otimes 1) \circ x^2$$
$$+ (j \otimes j + j \otimes k - 1 \otimes i - 1 \otimes 1 + i \otimes 1) \circ x$$
$$- k + j$$

*According to the theorem 8.5, we set* $q_2(x) = (x^2) \otimes 1$. *Then*

$$s_2(x) = r(x) - q_2(x) \circ (x - k) = r(x) - (x^2 \otimes 1) \circ (x - k)$$
$$= x^3 - x^2j - x^2k - jx^2 + jxj + jxk - xkx + xkj - x + ix - k + j$$
(10.1)
$$- x^3 + x^2k$$
$$= -x^2j - jx^2 + jxj + jxk - xkx + xkj - x + ix - k + j$$

*According to the theorem 8.5, we set*

(10.2) $$q_1(x) = -x \otimes j - jx \otimes 1 - xk \otimes 1$$

*The equality*

$$s_1(x) = s_2(x) - q_1(x) \circ (x - k)$$
$$= -x^2j - jx^2 + jxj + jxk - xkx + xkj - x + ix - k + j$$
$$- (-x \otimes j - jx \otimes 1 - xk \otimes 1) \circ (x - k)$$
(10.3)
$$= -x^2j - jx^2 + jxj + jxk - xkx + xkj - x + ix - k + j$$
$$+ x(x - k)j + jx(x - k) + xk(x - k)$$
$$= -x^2j - jx^2 + jxj + jxk - xkx + xkj - x + ix - k + j$$
$$+ x^2j - xkj + jx^2 - jxk + xkx - xkk$$
$$= jxj + ix - k + j$$

*follows from equalities* (10.1), (10.2). *According to the theorem 8.5, we set*

(10.4) $$q_0(x) = j \otimes j + i \otimes 1$$

*The equality*

$$s_0 = s_1(x) - q_0(x) \circ (x - k)$$
$$= jxj + ix - k + j - (j \otimes j + i \otimes 1) \circ (x - k)$$
(10.5)
$$= jxj + ix - k + j - j(x - k)j - i(x - k)$$
$$= jxj + ix - k + j - jxj + jkj - ix + ik$$
$$= -k + j + jkj + ik = 0$$

*follows from equalities* (10.3), (10.4). *The equality*

$$r(x) = s_0 + (q_0 + q_1(x) + q_2(x)) \circ (x - k)$$
(10.6)
$$= (j \otimes j + i \otimes 1 - x \otimes j - jx \otimes 1 - xk \otimes 1 + x^2 \otimes 1) \circ (x - k)$$

*follows from equalities* (10.1), (10.2), (10.4), (10.5). □

The quotient

(10.7) $$q(x) = j \otimes j + i \otimes 1 - x \otimes j - jx \otimes 1 - xk \otimes 1 + x^2 \otimes 1$$

of the polynomial

$$r(x) = (x - j)(x - k)(x - j - k)$$

divided by the polynomial

$$p(x) = x - k$$

is a tensor depending on $x$, not a polynomial. That is why we cannot use the theorem 8.5.

**Question 10.2.** *It is easy to see that*

$$q(j) = (j \otimes j + i \otimes 1 - x \otimes j - jx \otimes 1 - xk \otimes 1 + x^2 \otimes 1)(j)$$
$$= j \otimes j + i \otimes 1 - j \otimes j - j^2 \otimes 1 - jk \otimes 1 + j^2 \otimes 1 = 0 \quad (10.8)$$

*Therefore, the question arises of how to select in the tensor*

$$j \otimes j + i \otimes 1 - x \otimes j - jx \otimes 1 - xk \otimes 1 + x^2 \otimes 1 \quad (10.9)$$

*factor $x - j$.* □

**Example 10.3.** *One possible answer to the question 10.2 is the transformation*

$$j \otimes j + i \otimes 1 - x \otimes j - jx \otimes 1 - xk \otimes 1 + x^2 \otimes 1$$
$$= j \otimes j - x \otimes j + i \otimes 1 - jx \otimes 1 - xk \otimes 1 + x^2 \otimes 1 \quad (10.10)$$
$$= -(x - j) \otimes j + (i - jx - xk + x^2) \otimes 1$$

*Consider the polinomial*

$$i - jx - xk + x^2 \quad (10.11)$$

*According to the theorem 8.5, we set $q_1(x) = x \otimes 1$. Then*

$$s_1(x) = r(x) - q_1(x) \circ (x - j) = r(x) - (x \otimes 1) \circ (x - j)$$
$$= i - jx - xk + x^2 - x^2 + xj \quad (10.12)$$
$$= i - jx - xk + xj$$

*Now we set*

$$q_0 = -j \otimes 1 - 1 \otimes k + 1 \otimes j \quad (10.13)$$

*Then*

$$s_0 = s_1(x) - q_0 \circ p(x) = s_1(x) - (-j \otimes 1 - 1 \otimes k + 1 \otimes j) \circ (x - j)$$
$$= i - jx - xk + xj + j(x - j) + (x - j)k - (x - j)j$$
$$= i - jx - xk + xj + jx - j^2 + xk - jk - xj + j^2 \quad (10.14)$$
$$= i - jk = 0$$

*The equality*

$$r(x) = s_1(x) + (x \otimes 1) \circ (x - j)$$
$$= (-j \otimes 1 - 1 \otimes k + 1 \otimes j + x \otimes 1) \circ (x - j) \quad (10.15)$$

*follows from equalities* (10.12), (10.14). *The equality*

$$j \otimes j + i \otimes 1 - x \otimes j - jx \otimes 1 - xk \otimes 1 + x^2 \otimes 1$$
$$= -(x - j) \otimes j + ((-j \otimes 1 - 1 \otimes k + 1 \otimes j + x \otimes 1) \circ (x - j)) \otimes 1$$
$$= ((-1 \otimes j - j \otimes 1 - 1 \otimes k + 1 \otimes j + x \otimes 1) \circ (x - j)) \otimes 1 \quad (10.16)$$
$$= ((-j \otimes 1 - 1 \otimes k + x \otimes 1) \circ (x - j)) \otimes 1$$

*follows from equalities* (10.10), (10.15). □

We see that quotient in the example 10.3 is $H^{3\otimes}$-number.

**Example 10.4.** *The equality*

(10.17)
$$r(x) = (((-j \otimes 1 - 1 \otimes k + 1 \otimes j + x \otimes 1) \circ (x - j)) \otimes 1) \circ (x - k)$$
$$= (-j \otimes 1 \otimes 1 - 1 \otimes k \otimes 1 + 1 \otimes j \otimes 1 + x \otimes 1 \otimes 1) \circ (x - j, x - k)$$

*follows from equalities* (10.6), (10.16). □

**Question 10.5.** *The equality* (10.17) *is an answer to the question* 10.2. *We see that factorization of the polynomial* $r(x)$ *is bilinear map of polynomial* $x - j$, $x - k$. *However this answer is not complete. According to the remark* 2.6, *the representation* (8.5) *of polynomial is not unique. The question of the choice of factors and their order is open. At the same time there is relation between different representations of polynomial as product of factors.* □

If the polynomial $p(x)$ has form $p(x) = p_1 \circ x + p_0$ where $p_1$ be nonsingular tensor, then we can use the theorem 8.5 to divide by the polynomial $p(x)$, if we consider the equality

(10.18)
$$p(x) = p_1 \circ (x + p_1^{-1} \circ p_0)$$

**Question 10.6.** *The polynomial*

(10.19)
$$p_1(x) = ix - xi - 1$$

*is divisor of the polynomial*

(10.20)
$$\begin{aligned}
p(x) &= (x - j)(ix - xi - 1)(x - k) \\
&= ((x - j)ix - (x - j)xi - x + j)(x - k) \\
&= xix(x - k) - jix(x - k) - x^2 i(x - k) \\
&\quad + jxi(x - k) - x(x - k) + j(x - k) \\
&= xix^2 - xixk + kx^2 - kxk - x^2 ix + x^2 ik \\
&\quad + jxix - jxik - x^2 + xk + jx - jk \\
&= (1 \otimes i \otimes 1 \otimes 1 - 1 \otimes 1 \otimes i \otimes 1) \circ x^3 \\
&\quad + (j \otimes i \otimes 1 - 1 \otimes i \otimes k + k \otimes 1 \otimes 1 - 1 \otimes 1 \otimes j - 1 \otimes 1 \otimes 1) \circ x^2 \\
&\quad + (1 \otimes k + j \otimes (1 + j) \circ x + 1 - i
\end{aligned}$$

*However, we cannot use the theorem* 8.5 *to find quotient, because the tensor*

$$i \otimes 1 - 1 \otimes i \in H^{2\otimes}$$

*is singular. Therefore, we need additional research to solve this problem.* □

## 11. Division of Polynomials in Non-associative Algebra

**Remark 11.1.** *In non-associative algebra, brackets determine in what order we perform the multiplication. So monomial can be written as follows*

(11.1)
$$(p(x)q(x))r(x)$$

(11.2)
$$p(x)(q(x)r(x))$$

where $p$, $q$, $r$ are also monomials. This imposes a limitation on the possibility to divide polynomial by a polynomial. □

In the example 11.2, I considered the polynomial $r(x)$ of power 3 in octonion algebra. Becaus order of factors in product is important for defining structure of tensor, I will explicitly write brackets even when the order is evident. For instance, since, in general,

$$(11.3) \quad ((a\otimes)b) \circ x = (ax)b \ne a(xb) = (a(\otimes b)) \circ x$$

then expressions $(a\otimes)b$ and $a(\otimes b)$ represent different linear maps.

**Example 11.2.** *Let*
$$p(x) = x - k$$
$$\begin{aligned}r(x) &= ((x-j)(x-k))(x-jl) = (x^2 - jx - xk + i)(x - jl)\\
(11.4) \quad &= (x^2)(x-jl) - (jx)(x-jl) - (xk)(x-jl) + i(x-jl)\\
&= (x^2)x - (x^2)jl - (jx)x - (xk)x + (jx)jl + (xk)jl + ix - ijl\end{aligned}$$
*According to the theorem 8.5, we set* $q_2(x) = (x^2)(\otimes 1)$. *Then*
$$\begin{aligned}s_2(x) &= r(x) - q_2(x) \circ (x-k) = r(x) - ((x^2)(\otimes 1)) \circ (x-k)\\
&= (x^2)x - (x^2)jl - (jx)x - (xk)x + (jx)jl + (xk)jl + ix + kl\\
(11.5) \quad &\quad - (x^2)x + (x^2)k\\
&= (x^2)k - (x^2)jl - (jx)x - (xk)x + (jx)jl + (xk)jl + ix + kl\end{aligned}$$
*According to the theorem 8.5, we set*
$$(11.6) \quad q_1(x) = (x\otimes)k - (x\otimes)jl - (jx)(\otimes 1) - (xk)(\otimes 1)$$
*The equality*
$$\begin{aligned}s_1(x) &= s_2(x) - q_1(x) \circ (x-k)\\
&= (x^2)k - (x^2)jl - (jx)x - (xk)x + (jx)jl + (xk)jl + ix + kl\\
&\quad - ((x\otimes)k - (x\otimes)jl - (jx)(\otimes 1) - (xk)(\otimes 1)) \circ (x-k)\\
(11.7) \quad &= (x^2)k - (x^2)jl - (jx)x - (xk)x + (jx)jl + (xk)jl + ix + kl\\
&\quad - (x(x-k))k + (x(x-k))jl + (jx)(x-k) + (xk)(x-k)\\
&= (x^2)k - (x^2)jl - (jx)x - (xk)x + (jx)jl + (xk)jl + ix + kl\\
&\quad - (x^2)k + (xk)k + (x^2)jl - (xk)jl + (jx)x - (jx)k + (xk)x - (xk)k\\
&= (jx)(jl-k) + ix + kl\end{aligned}$$
*follows from equalities* (11.5), (11.6). *According to the theorem 8.5, we set*
$$(11.8) \quad q_0(x) = (j\otimes)(jl-k) + i(\otimes 1)$$
*The equality*
$$\begin{aligned}s_0 &= s_1(x) - q_0(x) \circ (x-k)\\
&= (jx)(jl-k) + ix + kl - ((j\otimes(jl-k)) + i(\otimes 1)) \circ (x-k)\\
(11.9) \quad &= (jx)(jl-k) + ix + kl - (j(x-k))(jl-k) - i(x-k)\\
&= (jx)(jl-k) + ix + kl - (jx - jk)(jl-k) - ix + ik\\
&= (jx)(jl-k) + ix + kl - (jx)(jl-k) + i(jl-k) - ix - j\\
&= kl + ijl - ik - j = 0\end{aligned}$$

*follows from equalities* (11.7), (11.8). *The equality*

$$r(x) = s_0 + (q_0 + q_1(x) + q_2(x)) \circ (x - k)$$

(11.10)
$$= ((j\otimes)(jl - k) + i(\otimes 1)$$
$$+ (x\otimes)k - (x\otimes)jl - (jx)(\otimes 1) - (xk)(\otimes 1) + (x^2)(\otimes 1)) \circ (x - k)$$

*follows from equalities* (11.5), (11.6), (11.8), (11.9). □

**Example 11.3.** *The quotient of the polynomial* (11.4) *divided by the polynomial*

$$p(x) = x - k$$

*has the following form*

(11.11)
$$q(x) = (j\otimes)(jl - k) + i(\otimes 1)$$
$$+ (x\otimes)k - (x\otimes)jl - (jx)(\otimes 1) - (xk)(\otimes 1) + (x^2)(\otimes 1)$$
$$= ((x - j)\otimes)(k - jl) + (x^2 - jx - xk + i)(\otimes 1)$$

(11.4) $\boxed{r(x) = ((x - j)(x - k))(x - jl)}$

*The first term of the tensor* $q(x)$ *has factor* $x - j$. *To check if the second term has factor* $x - j$, *consider the polynomial*

(11.12) $$r_1(x) = x^2 - jx - xk + i$$

*According to the theorem* 8.5, *we set* $q_1(x) = x(\otimes 1)$. *Then*

(11.13)
$$s_1(x) = r_1(x) - (x(\otimes 1)) \circ (x - j)$$
$$= x^2 - jx - xk + i - x(x - j) = x^2 - jx - xk + i - x^2 + xj$$
$$= -jx - xk + xj + i$$

*Now we set*

(11.14) $$q_0 = -j(\otimes 1) - (1\otimes)k + (1\otimes)j$$

*Then*

(11.15)
$$s_0 = s_1(x) - q_0 \circ p(x) = s_1(x) - (-j(\otimes 1) - (1\otimes)k + (1\otimes)j) \circ (x - j)$$
$$= -jx - xk + xj + i + j(x - j) + (x - j)k - (x - j)j$$
$$= -jx - xk + xj + i + jx + 1 + xk - i - xj - 1 = 0$$

*The equality*

(11.16)
$$r_1(x) = s_1(x) + (x(\otimes 1)) \circ (x - j)$$
$$= (-j(\otimes 1) - (1\otimes)k + (1\otimes)j + x(\otimes 1)) \circ (x - j)$$
$$= ((1\otimes)(j - k) + (x - j)(\otimes 1)) \circ (x - j)$$

*follows from equalities* (11.13), (11.15). □

From equalities (11.11), (11.16), it follows that the quotient of the polynomial (11.4) divided by the polynomial

$$p(x) = x - k$$

has the following form

$$\begin{aligned} q(x) = &= ((x-j)\otimes)(k-jl) \\ &+ (((1\otimes)(j-k) + (x-j)(\otimes 1)) \circ (x-j))(\otimes 1) \\ &= (((1\otimes_1)1\otimes)(k-jl) \\ &+ ((1\otimes_1)(j-k) + (x-j)(\otimes_1 1))(\otimes 1)) \circ (x-j) \\ &= (((1\otimes_1)1\otimes)(k-jl) \\ &+ ((1\otimes_1)(j-k))(\otimes 1) + ((x-j)(\otimes_1 1))(\otimes 1)) \\ &\circ (x-j) \end{aligned} \quad (11.17)$$

(11.4) $\boxed{r(x) = ((x-j)(x-k))(x-jl)}$

Therefore, we can represent the polynomial (11.4) as follows

$$\begin{aligned} (11.18) \quad r(x) = &(((1\otimes_1)1\otimes_2)(k-jl) \\ &+ ((1\otimes_1)(j-k))(\otimes_2 1) + ((x-j)(\otimes_1 1))(\otimes_2 1)) \\ &\circ (x-j, x-k) \end{aligned}$$

## 13. Index

# 14. Special Symbols and Notations

# Многочлен в некомутативной алгебре

Александр Клейн

**Аннотация.** Я рассмотрел определение и свойства многочлена в некомутативной алгебре. Многочлен может иметь конечное или бесконечное, а также пустое множество корней. Например, многочлен

$$p_1(x) = ix - xi - 1$$

не имеет корней и многочлен

$$p_k(x) = ix - xi - k$$

имеет множество корней

$$x = C_1 + C_2 i + \frac{1}{2} j$$

Рассмотрена операция деления многочленов с остатком.

## Содержание

## 1. Линейное отображение

**Соглашение 1.1.** *Если отображение*

$$f : A \to A$$

Aleks_Kleyn@MailAPS.org.
  http://AleksKleyn.dyndns-home.com:4080/   http://arxiv.org/a/kleyn_a_1.
  http://AleksKleyn.blogspot.com/.

$D$-алгебры $A$ является линейным отображением, то я пользуюсь обозначением

$$(1.1) \qquad f \circ a = f(a)$$

для образа отображения $f$. □

**Определение 1.2.** *Линейное отображение*
$$f : A \to A$$
*$D$-алгебры $A$ удовлетворяет равенствам*

$$(1.2) \qquad f \circ (a + b) = f \circ a + f \circ b$$

$$(1.3) \qquad f \circ (da) = d(f \circ a)$$

$$a, b \in A \quad d \in D$$

*Обозначим* $\mathcal{L}(D; A \to A)$ *множество линейных отображений $D$-алгебры $A$.*
□

Пусть $H$ - алгебра кватернионов.

**Теорема 1.3.** *Пусть произведение в алгебре $H \otimes H$ определено согласно правилу*

$$(1.4) \qquad (p_0 \otimes p_1) \circ (q_0 \otimes q_1) = (p_0 q_0) \otimes (q_1 p_1)$$

*Представление*

$$(1.5) \qquad h : H \otimes H \longrightarrow\!\!\!* \mathcal{L}(R; H \to H) \quad h(p) : g \to p \circ g$$

*$R$-алгебры $H \otimes H$ в абелевой группе $\mathcal{L}(R; H \to H)$, определённое равенством*

$$(1.6) \qquad (a \otimes b) \circ g = agb \quad a, b \in H \quad g \in \mathcal{L}(R; H \to H)$$

*является левым $H \otimes H$-модулем. Если мы положим $g = E$, где $E \in \mathcal{L}(R; H \to H)$ - тождественное отображение, то мы можем отождествить линейное отображение $f : H \to H$ алгебры кватернионов и тензор*

$$(1.7) \qquad f = f_{s \cdot 0} \otimes f_{s \cdot 1} \in H \otimes H$$

*опираясь на равенство*

$$(1.8) \qquad f \circ a = (f_{s \cdot 0} \otimes f_{s \cdot 1}) \circ a = f_{s \cdot 0}\, a\, f_{s \cdot 1}$$

Доказательство. Теорема является следствием утверждения что отображение
$$E \circ x = x$$
порождает любое линейное отображение алгебры кватернионов □

**Определение 1.4.** *Выражение* $f_{s \cdot p}$, $p = 0, 1$, *в равенстве* $(1.7)$ *называется* **компонентой линейного отображения** $f$. □

$$\boxed{(1.7) \qquad f = f_{s \cdot 0} \otimes f_{s \cdot 1} \in H \otimes H}$$

**Определение 1.5.** *Полилинейное отображение D-алгебры A*
$$f: A^n \to A$$
*удовлетворяет равенствам*

(1.9) $\quad f \circ (a_1, ..., a_i + b_i, ..., a_n) = f \circ (a_1, ..., a_i, ..., a_n) + f \circ (a_1, ..., b_i, ..., a_n)$

(1.10) $\quad\quad\quad f \circ (a_1, ..., pa_i, ..., a_n) = pf \circ (a_1, ..., a_i, ..., a_n)$

$$1 \le i \le n \quad a_i, b_i \in A \quad p \in D$$

*Обозначим* $\mathcal{L}(D; A^n \to A)$ *множество n-линейных отображений D-алгебры A.* □

**Теорема 1.6.** *Представление*
$$h: H^{n+1\otimes} \times S_n \longrightarrow \mathcal{LA}(R; H^n \to H)$$
*алгебры* $H^{n+1\otimes}$ *в модуле* $\mathcal{LA}(R; H^n \to H)$, *определённое равенством*

(1.11) $\quad (a_0 \otimes ... \otimes a_n, \sigma) \circ (f_1 \otimes ... \otimes f_n) = a_0\sigma(f_1)a_1...a_{n-1}\sigma(f_n)a_n$

$$a_0, ..., a_n \in H \quad \sigma \in S_n \quad f_1, ..., f_n \in \mathcal{L}(R; H \to H)$$

*позволяет отождествить тензор* $d \in A^{n+1\otimes}$ *и перестановку* $\sigma \in S^n$ *с отображением*

(1.12) $\quad\quad\quad (d, \sigma) \circ (f_1, ..., f_n) \quad f_i = \delta \in \mathcal{L}(R; H \to H)$

*где* $\delta \in \mathcal{L}(R; H \to H)$ - *тождественное отображение.*

## 2. Многочлен

**Определение 2.1.** *Пусть*
$$f: A^n \to A$$
*полилинейное отображение D-модуля A.e Отображение*
$$p: A \to A$$
*определённое равенством*

(2.1) $\quad\quad\quad p \circ x^n = f \circ (x_1, ..., x_n) \quad x_1 = ... = x_n = x$

*называется* **однородным многочленом** *степени n.* □

**Теорема 2.2.** *Множество* $A_n[x]$ *однородных многочленов степени n является абелевой группой.*

Доказательство. Теорема является следствием утверждения, что сумма $n$-линейных отображений является $n$-линейным отображением. □

**Теорема 2.3.** *Пусть A - D-алгебра. Тогда для любого однородного многочлена*
$$p: A \to A$$

*степени $n$ существует тензор $p' \in A^{n+1\otimes}$ такой, что*

$$(2.2) \quad p \circ x^n = p' \circ x^n = p' \circ (x_1, ..., x_n) \quad x_1 = ... = x_n = x$$

Доказательство. Теорема является следствием теоремы 1.6. В дальнейшем мы будем отождествлять однородный многочлен $p$ и соответствующий тензор $p'$.  □

**Определение 2.4.** *Если тензор $p$ имеет вид $p = p_0 \otimes ... \otimes p_n$, то многочлен $p \circ x^n$ называется одночленом.*  □

**Теорема 2.5.** *Пусть $p_k(x)$ - **одночлен степени** $k$ над ассоциативной $D$-алгеброй $A$. Тогда*

2.5.1: *Одночлен степени $0$ имеет вид $p_0(x) = a_0$, $a_0 \in A$.*

2.5.2: *Если $k > 0$, то*

$$(2.3) \quad p_k(x) = p_{k-1}(x) x a_k$$

*где $a_k \in A$.*

Доказательство. Мы докажем утверждение теоремы индукцией по степени $n$ одночлена.

Пусть $n = 0$. Так как одночлен $p_0(x)$ является константой, то мы получаем утверждение 2.5.1.

Пусть $n = k$. Последний множитель одночлена $p_k(x)$ является либо $a_k \in A$, либо имеет вид $x^l$, $l \geq 1$. В последнем случае мы положим $a_k = 1$. Множитель, предшествующий $a_k$, имеет вид $x^l$, $l \geq 1$. Мы можем представить этот множитель в виде $x^{l-1}x$. Следовательно, утверждение доказано.  □

**Замечание 2.6.** *В теореме 2.5 я рассмотрел рекурсивное представление одночлена в ассоциативной $D$-алгебре. Так как скобки в произведении я могу расставить произвольно, то в ассоциативной $D$-алгебре рекурсивное представление одночлена неоднозначно. Например, одночлен $ax^2bxcx^3d$ может быть записан в любой из следующих форм*

$$(2.4) \quad \begin{aligned} ax^2bxcx^3d &= ax(xbxcx^3) = (ax)x(bxcx^3d) = (ax^2b)x(cx^3d) \\ &= (ax^2bxc)x(x^2d) = (ax^2bxcx)x(xd) = (ax^2bxcx^2)xd \end{aligned}$$

*Я выбрал равенство (2.3) как самое простое для алгоритма деления многочленов.*  □

**Определение 2.7.** *Обозначим*

$$A[x] = \bigoplus_{n=0}^{\infty} A_n[x]$$

*прямую сумму[2.1] абелевых групп $A_n[x]$. Элемент $p(x)$ абелевой группы $A[x]$ называется **многочленом** над $D$-алгеброй $A$.*  □

---

[2.1] Смотри определение прямой суммы абелевых групп в [1], страница 55. Согласно теореме 1 на той же странице, прямая сумма абелевых групп существует.

Следовательно, многочлен степени $n$ может быть записан в виде

$$(2.5) \qquad p(x) = p_0 + p_1 \circ x + ... + p_n \circ x^n \quad p_i \in A^{i+1\otimes} \quad i = 0,...,n$$

## 3. Операции с многочленами

**Определение 3.1.** *Пусть*

$$p(x) = p_0 + p_1 \circ x + ... + p_n \circ x^n \quad p_i \in A^{i+1\otimes} \quad i = 0,...,n$$

$$r(x) = r_0 + r_1 \circ x + ... + r_n \circ x^n \quad r_i \in A^{i+1\otimes} \quad i = 0,...,n$$

*многочлены.*[3.1] *Мы определим сумму многочленов $p$ и $r$ равенством*

$$(3.1) \qquad (p+r)(x) = p_0 + r_0 + (p_1 + r_1) \circ x + ... + (p_n + r_n) \circ x^n$$

□

**Определение 3.2.** *Билинейное отображение*

$$(3.2) \qquad * : A^{n\otimes} \times A^{m\otimes} \to A^{n+m-1\otimes}$$

*определено равенством*

$$(3.3) \qquad (a_1 \otimes ... \otimes a_n) * (b_1 \otimes ... \otimes b_n) = a_1 \otimes ... \otimes a_{n-1} \otimes a_n b_1 \otimes b_2 \otimes ... \otimes b_n$$

□

**Теорема 3.3.** *Для любых тензоров $a \in A^{n+1\otimes}$, $b \in A^{m+1\otimes}$, произведение однородных многочленов $a \circ x^n$, $b \circ x^m$ определено равенством*

$$(3.4) \qquad (a*b) \circ x^{n+m} = (a \circ x^n)(b \circ x^m)$$

## 4. Линейное уравнение

Пусть $\overline{\overline{e}}$ - базис конечно мерной алгебры $A$ над полем $F$ и $C_{ij}^k$ - структурные константы алгебры $A$ относительно базиса $\overline{\overline{e}}$.

Рассмотрим линейное уравнение

$$(4.1) \qquad a \circ x = b$$

где $a = a_{s \cdot 0} \otimes a_{s \cdot 1} \in A^{2\otimes}$. Согласно теореме [3]-6.4.1, мы можем записать уравнение (4.1) в стандартной форме

$$(4.2) \qquad a^{ij} e_i x e_j = b$$

$$(4.3) \qquad a^{ij} = a_{s \cdot 0}{}^i a_{s \cdot 1}{}^j \quad a_{s \cdot 0} = a_{s \cdot 0}{}^i e_i \quad a_{s \cdot 1} = a_{s \cdot 1}{}^i e_i$$

Согласно теореме [3]-6.4.5, уравнение (4.2) эквивалентно уравнению

$$(4.4) \qquad a_i^j x^i = b^j$$

$$(4.5) \qquad a_i^j = a^{kr} C_{ki}^p C_{pr}^j \quad x = x^i e_i \quad b = b^i e_i$$

---

[3.1] Если коэффициент, допустим $p_i$, отсутствует, то мы полагаем $p_i = 0$.

Согласно теории линейных уравнений над полем, если определитель

(4.6) $$\det \|a_i^j\| \neq 0$$

то уравнение (4.1) имеет единственное решение.

(4.1)    $a \circ x = b$

[3] Александр Клейн, Линейное отображение $D$-алгебры, eprint arXiv:1502.04063 (2015)

**Определение 4.1.** *Тензор $a \in A^{2\otimes}$ называется* **невырожденным тензором**, *если этот тензор удовлетворяет условию* (4.6). □

**Теорема 4.2.** *Пусть $a \in A^{2\otimes}$ - невырожденный тензор. Если равенство (4.2) рассматривать как преобразование алгебры $A$, то обратное преобразование можно записать в виде*

(4.7) $$x = c^{pq} e_p b e_q$$

*где компоненты $c^{pq}$ удовлетворяют уравнению*

(4.8) $$\delta_0^r \delta_0^s = a^{ij} c^{pq} C_{ip}^r C_{qj}^s$$

Доказательство. Теорема является следствием теоремы [2]-6.7. □

(4.1)    $a \circ x = b$     (4.2)    $a^{ij} e_i x e_j = b$

[2] Александр Клейн, Многочлен над ассоциативной $D$-алгеброй, eprint arXiv:1302.7204 (2013)

**Определение 4.3.** *Пусть $a \in A^{2\otimes}$ - невырожденный тензор. Тензор*

(4.9) $$a^{-1} = c^{pq} e_p \otimes e_q$$

*называется* **тензором, обратным тензору** $a$. □

**Теорема 4.4.** *Если $a \in A^{2\otimes}$ - невырожденный тензор, то линейное уравнение*

(4.10) $$a \circ x = b$$

*имеет единственное решение*

(4.11) $$x = a^{-1} \circ b$$

*Если $a$- невырожденный тензор, то уравнение (4.10) имеет решение только при условии $b \in \ker a$. В этом случае уравнение (4.10) имеет бесконечно много решений.*

(4.10)    $a \circ x = b$

**Пример 4.5.** *Тензор*
$$i \otimes 1 - 1 \otimes i \in H^{2\otimes}$$

*соответствует матрице*

$$(4.12) \quad \begin{pmatrix} 0 & -1 & 0 & 0 \\ 1 & 0 & 0 & 0 \\ 0 & 0 & 0 & -1 \\ 0 & 0 & 1 & 0 \end{pmatrix} - \begin{pmatrix} 0 & -1 & 0 & 0 \\ 1 & 0 & 0 & 0 \\ 0 & 0 & 0 & 1 \\ 0 & 0 & -1 & 0 \end{pmatrix} = \begin{pmatrix} 0 & 0 & 0 & 0 \\ 0 & 0 & 0 & 0 \\ 0 & 0 & 0 & -2 \\ 0 & 0 & 2 & 0 \end{pmatrix}$$

*Следовательно многочлен*

$$(4.13) \quad p_1(x) = ix - xi - 1$$

*не имеет корней и многочлен*

$$(4.14) \quad p_k(x) = ix - xi - k$$

*имеет множество корней*

$$(4.15) \quad x = C_1 + C_2 i + \frac{1}{2} j$$

□

## 5. Левостороний многочлен

[6] Paul M. Cohn, Skew Fields, Cambridge University Press, 1995

В книге [6] на странице 48, профессор Кон писал, что изучать многочлены над не коммутативным полем - сложная задача. Поэтому профессор Оре рассмотрел кольцо правых косых многочленов[5.1] где Оре определил преобразование одночлена

$$(5.1) \quad at = ta^\alpha + a^\delta$$

Отображение

$$\alpha : a \in A \to a^\alpha \in A$$

является эндоморфизмом $D$-алгебры $A$ и отображение

$$\delta : a \in A \to a^\delta \in A$$

является $\alpha$-дифференцированием $D$-алгебры $A$

$$(5.2) \quad (a+b)^\delta = a^\delta + b^\delta \quad (ab)^\delta = a^\delta b^\alpha + ab^\delta$$

В этом случае все коэффициенты многочлена могут быть записаны справа. Подобные многочлены называются правосторонними многочленами.

Аналогично мы можем рассматривать преобразование одночлена

$$(5.3) \quad ta = a^\alpha t + a^\delta$$

В этом случае все коэффициенты многочлена могут быть записаны слева. Подобные многочлены называются левосторонними многочленами.[5.2]

Точка зрения, предложенная Оре, позволила получить качественную картину теории полиномов в некоммутативной алгебре. Однако существуют алгебры где рассматриваемая картина не полна.

---

[5.1] Смотри также определение на страницах [6]-(48 - 49).

[5.2] Смотри также определение на странице [5]-3.

**Пример 5.1.** *Если преобразование* (5.1) *определено в D-алгебре A, то мы можем любой многочлен p степени* 1 *отобразить в правосторонний многочлен p′ степени* 1 *таким образом, что*

$$p(t) = p'(t)$$

*Многочлен p′ может иметь следующий вид*

(5.4) $$p'(t) = ta_1 + a_0 \quad a_1 \neq 0$$

(5.5) $$p'(t) = t\,0 + a_0 \quad a_0 \neq 0$$

(5.6) $$p'(t) = t\,0 + 0$$

*Многочлен* (5.4) *имеет единственное решение. Многочлен* (5.5) *не имеет решений. Любое A-число является корнем многочлена* (5.6).

*Очевидно, что мы упустили случай, когда многочлен p имеет бесконечно много корней, но не любое A-число является корнем многочлена p. Поэтому возникает вопрос о существовании преобразования* (5.1) *в алгебре кватернионов.* □

**Пример 5.2.** *Пусть*

$$p(x) = \sum_{i=0}^{n} p_i x^i$$

*и*

$$r(x) = \sum_{j=0}^{m} r_j x^j$$

*левосторонние многочлены. Мы определим произведение многочленов p и r равенством*[5.3]

(5.7) $$p(x) * r(x) = (p * r)(x) = \sum_{j=0}^{m+n} (\sum_{i=0}^{j} p_i r_{j-i}) x^j$$

*Рассмотрим многочлены*

(5.8) $$L(x) = x - i$$

(5.9) $$R(x) = x - j$$

*Равенство*

(5.10) $$P(x) = (L * R)(x) = x^2 - (i+j)x + k$$

*является следствием равенств* (5.7), (5.8), (5.9). *Пусть*

$$h = R(i) = i - j$$

$$\tilde{x} = hxh^{-1} = (i-j)x(i-j)^{-1} = \frac{1}{2}(i-j)x(j-i)$$

(5.11) $$L(\tilde{x}) = \frac{1}{2}(i-j)x(j-i) - i$$

*Согласно теореме* [5]*-1, пункт* (ii), *многочлен*

(5.12) $$P_1(x) = L(\tilde{x})R(x) = (\frac{1}{2}(i-j)x(j-i) - i)(x - j)$$

*совпадает с многочленом* $P(x)$, $P_1(x) = P(x)$

(5.13) $$x^2 - (i+j)x + k = (\frac{1}{2}(i-j)x(j-i) - i)(x - j)$$

*Однако равенство* (5.13) *не верно. Пусть* $x = i + j$. *Равенства*

(5.14) $$(i+j)^2 - (i+j)(i+j) + k = (\frac{1}{2}(i-j)(i+j)(j-i) - i)(i + j - j)$$

(5.15) $$\begin{aligned} k &= (\frac{1}{2}(i^2 + ij - ji - j^2)(j - i) - i)i \\ &= (k(j-i) - i)i = (-i - j - i)i \end{aligned}$$

*являются следствием равенства* (5.13). □

## 6. Квадратный корень

**Определение 6.1.** *Решение* $x = \sqrt{a}$ *уравнения*

(6.1) $$x^2 = a$$

*в D-алгебре A называется* **квадратным корнем** *A-числа* $a$. □

**Теорема 6.2.** *Пусть H - алгебра кватернионов и a - H-число.*
  6.2.1: *Если* $\operatorname{Re} \sqrt{a} \neq 0$, *то уравнение*
$$x^2 = a$$
  *имеет корни* $x = x_1$, $x = x_2$ *такие, что*

(6.2) $$x_2 = -x_1$$

  6.2.2: *Если* $a = 0$, *то уравнение*
$$x^2 = a$$
  *имеет корень* $x = 0$ *кратности* 2.
  6.2.3: *Если условия* 6.2.1, 6.2.2 *не верны, то уравнение*
$$x^2 = a$$
  *имеет бесконечно много корней таких, что*

(6.3) $$x \in \operatorname{Im} H \quad a \in \operatorname{Re} H \quad |x| = \sqrt{-a}$$

Доказательство. Теорема является следствием теоремы [4]-8.5. □

**Теорема 6.3.** *Пусть кватернион a является корнем многочлена*

(6.4) $$r(x) = x^2 + 1$$

*Тогда*

(6.5) $$\sqrt{a} = \pm \frac{1}{\sqrt{2}}(1 + a)$$

---

[5.3] Смотри также определение на странице [5]-3.

Доказательство. Равенство

$$(6.6) \qquad (a+1)^2 = a^2 + 2a + 1 = 2a$$

следует из утверждения, что $a$ является корнем многочлена (6.4). Равенство (6.5) является следствием равенства (6.6) и утверждения 6.2.1. □

## 7. Многочлен с заданными корнями

В коммутативной алгебре, если я имею два корня многочлена второй степени, то этот многочлен определён однозначно. Ситуация несколько иная в некоммутативной алгебре.

**Теорема 7.1.** *Пусть* $x_1, x_2 \in A$,

$$(7.1) \qquad x_1 x_2 \neq x_2 x_1$$

*Тогда*

$$(7.2) \qquad p_{12}(x) = (x - x_1)(x - x_2) = x^2 - x_1 x - x x_2 + x_1 x_2$$

$$(7.3) \qquad p_{21}(x) = (x - x_2)(x - x_1) = x^2 - x_2 x - x x_1 + x_2 x_1$$

*неравные многочлены*

$$(7.4) \qquad p_{12} \neq p_{21}$$

Доказательство. Пусть

$$(7.5) \qquad x = x_1 + x_2$$

Равенство

$$(7.6) \qquad p_{12}(x_1 + x_2) = (x_1 + x_2 - x_1)(x_1 + x_2 - x_2) = x_2 x_1$$

является следствием равенств (7.2), (7.5). Равенство

$$(7.7) \qquad p_{21}(x_1 + x_2) = (x_1 + x_2 - x_2)(x_1 + x_2 - x_1) = x_1 x_2$$

является следствием равенств (7.3), (7.5). Утверждение (7.4) является следствием равенств (7.6), (7.7) и утверждения (7.1). □

**Теорема 7.2.** *Пусть* $a_1, a_2 \in A \otimes A$. *Тогда многочлен*

$$(7.8) \qquad p(x) = a_1 \circ p_{12}(x) + a_2 \circ p_{21}(x)$$

*имеет корни* $x = x_1$, $x = x_2$.

Доказательство. Теорема является следствием утверждения, что многочлены $p_{12}$, $p_{21}$ имеют корни $x = x_1$, $x = x_2$. □

**Вопрос 7.3.** *Имеет ли многочлен (7.8) корни отличные от значений* $x = x_1$, $x = x_2$? □

**Вопрос 7.4.** *Существуют ли другие многочлены, которые имеют корни* $x = x_1$, $x = x_2$? □

Похоже, что ответ на вопрос 7.4 положителен.

**Пример 7.5.** *Рассмотрим множество многочленов*

$$(7.9) \qquad p(x) = a_1 \circ ((x-i)(x-j)) + a_2 \circ ((x-j)(x-i))$$

*которые, согласно теореме 7.2, имеют корни $x = i$, $x = j$. Многочлен*

$$(7.10) \qquad r(x) = x^2 + 1$$

*также имеет корни $x = i$, $x = j$.* □

**Вопрос 7.6.** *Принадлежит ли многочлен (7.10) множеству многочленов (7.9)?* □

Теорема 7.7 отвечает на вопрос 7.6.

**Теорема 7.7.** *Не существует тензоров $a_1, a_2 \in A \otimes A$ таких, что*

$$(7.11) \qquad a_1 \circ (x^2 - ix - xj + k) + a_2 \circ (x^2 - jx - xi - k) = x^2 + 1$$

Доказательство. Предположим, что теорема не верна.

**Утверждение 7.8.** *Пусть существуют тензоры $a_1, a_2 \in A \otimes A$ таких, что равенство (7.11) верно.* ⊙

Равенства

$$(7.12) \qquad a_1 \circ x^2 + a_2 \circ x^2 = x^2$$

$$(7.13) \qquad a_1 \circ (ix) + a_1 \circ (xj) + a_2 \circ (jx) + a_2 \circ (xi) = 0$$

$$(7.14) \qquad a_1 \circ k - a_2 \circ k = 1$$

следуют из равенства (7.11).

Пусть $x = 1$. Равенство

$$(7.15) \qquad a_1 \circ i + a_1 \circ j + a_2 \circ j + a_2 \circ i = 0$$

следуют из равенств (7.13).

Согласно теореме 6.3, существует $x$ такой, что $x^2 = i$. Равенство

$$(7.16) \qquad a_1 \circ i + a_2 \circ i = i$$

следует из равенства (7.12).

Согласно теореме 6.3, существует $x$ такой, что $x^2 = j$. Равенство

$$(7.17) \qquad a_1 \circ j + a_2 \circ j = j$$

следует из равенства (7.12).

Равенство

$$(7.18) \qquad a_1 \circ i + a_2 \circ i + a_1 \circ j + a_2 \circ j = i + j$$

следует из равенств (7.16), (7.17) и противоречит равенству (7.15). Следовательно, утверждение 7.8 не верно и теорема доказана. □

8. Деление с остатком

**Определение 8.1.** *$A$-число $a$ называется **левым делителем** $A$-числа $b$, если существует $A$-число $c$ такое, что*

(8.1) $$ac = b$$

□

**Определение 8.2.** *$A$-число $a$ называется **правым делителем** $A$-числа $b$, если существует $A$-число $c$ такое, что*

(8.2) $$ca = b$$

□

Симметрия между определениями 8.1 и 8.2 очевидна. Также как очевидно различие между левым и правым делителями в связи с некоммутативностью произведения. Однако мы можем рассмотреть определение, обобщающее определения 8.1 и 8.2.

**Определение 8.3.** *$A$-число $a$ называется **делителем** $A$-числа $b$, если существует $A \otimes A$-число $c$ такое, что*

(8.3) $$c \circ a = b$$

*$A \otimes A$-число $c$ называется **частным от деления** $A$-числа $b$ на $A$-число $a$.* □

**Определение 8.4.** *Пусть деление в $D$-алгебре $A$ не всегда определено. $A$-**число** $a$ **делит** $A$-**число** $b$ **с остатком**, если следующее равенство верно*

(8.4) $$c \circ a + f = b$$

*$A \otimes A$-число $c$ называется **частным от деления** $A$-числа $b$ на $A$-число $a$. $A$-число $d$ называется **остатком от деления** $A$-числа $b$ на $A$-число $a$.* □

**Теорема 8.5.** *Пусть*
$$p(x) = x - a$$
*- многочлен степени 1. Пусть*
$$r(x) = r_0 + r_1 \circ x + ... + r_k \circ x^k$$
*многочлен степени $k > 0$. Тогда*

(8.5) $$r(x) = s_0 + (q_0 + q_1(x) + ... + q_{k-1}(x)) \circ p(x)$$

*где $q_i \in A_i[x] \otimes A$, $i = 0, ..., k-1$, - однородный многочлен степени $i$.*

Доказательство. Пусть

(8.6) $$r_k = r_{k.0.s} \otimes ... \otimes r_{k.k.s}$$

Согласно теореме 2.5,

(8.7) $$\begin{aligned}r_k \circ x^k &= ((r_{k.0.s} \otimes ... \otimes r_{k.k-1.s}) \circ x^{k-1}) x r_{k.k.s} \\ &= (((r_{k.0.s} \otimes ... \otimes r_{k.k-1.s}) \circ x^{k-1}) \otimes r_{k.k.s}) \circ x\end{aligned}$$

Положим

(8.8) $$q_{k-1}(x) = ((r_{k.0.s} \otimes ... \otimes r_{k.k-1.s}) \circ x^{k-1}) \otimes r_{k.k.s}$$

Для доказательства теоремы достаточно обратить внимание, что степень многочлена

(8.9) $$s_{k-1}(x) = r(x) - q_{k-1}(x) \circ (x - a)$$

меньше $k$. Равенство

(8.10) $$r(x) = s_{k-1}(x) + q_{k-1}(x) \circ (x - a)$$

является следствием равенства (8.9). Мы применим рассуждения доказательства к многочлену $s_{k-1}(x)$. После конечного числа шагов мы получим многочлен $s_0$ степени 0. Равенство (8.5) является следствием равенства (8.10). □

Алгоритм, рассмотренный в доказательстве теоремы называется стандартный алгоритм деления многочленов.

## 9. Примеры деления многочленов

В последующих примерах мы применим стандартный алгоритм деления многочленов.

**Пример 9.1.** *Пусть*

$$p(x) = x - i$$

(9.1) $$r(x) = (x - j)(x - i) = x^2 - jx - xi - k$$

*Согласно теореме* 8.5, *положим* $q_1(x) = x \otimes 1$. *Тогда*

(9.2) $$\begin{aligned} s_1(x) &= r(x) - q_1(x) \circ p(x) = r(x) - (x \otimes 1) \circ (x - i) \\ &= x^2 - jx - xi - k - x(x - i) = x^2 - jx - xi - k - x^2 + xi \\ &= -jx - k \end{aligned}$$

*Теперь мы положим* $q_0 = -j \otimes 1$. *Тогда*

(9.3) $$\begin{aligned} s_0 &= s_1(x) - q_0 \circ p(x) = s_1(x) - (-i \otimes 1 - 1 \otimes j + 1 \otimes i) \circ (x - i) \\ &= -ix - xj + k + xi + i(x - i) + (x - i)j - (x - i)i \\ &= -ix - xj + k + xi + ix + 1 + xj - k - xi - 1 = 0 \end{aligned}$$

*Равенство*

(9.4) $$\begin{aligned} r(x) &= s_1(x) + (x \otimes 1) \circ (x - i) \\ &= (-j \otimes 1 + x \otimes 1) \circ (x - i) \end{aligned}$$

*является следствием равенств* (9.2), (9.3). *Мы можем упростить выражение в равенстве* (9.4)

(9.5) $$r(x) = ((x - j) \otimes 1) \circ (x - i) = (x - j)(x - i)$$

*Мы видим, что представление многочлена* $r(x)$ *в равенстве* (9.5) *совпадает с представление многочлена* $r(x)$ *в равенстве* (9.1). □

Ответ (9.5) в примере 9.1 очевиден, так как я намеренно выбрал делитель $p(x)$ равным второму множителю многочлена $r(x)$.

$\boxed{p(x) = x - i}$ $\boxed{r(x) = (x - j)(x - i) = x^2 - jx - xi - k}$

Пример 9.2 более интересен, так как я изменил порядок множителей.

**Пример 9.2.** *Пусть*
$$p(x) = x - i$$
(9.6) $$r(x) = (x-i)(x-j) = x^2 - ix - xj + k$$

*Согласно теореме* 8.5, *положим* $q_1(x) = x \otimes 1$. *Тогда*

(9.7) $$\begin{aligned}s_1(x) &= r(x) - (x \otimes 1) \circ (x - i) \\ &= x^2 - ix - xj + k - x(x-i) = x^2 - ix - xj + k - x^2 + xi \\ &= -ix - xj + k + xi\end{aligned}$$

*Теперь мы положим*

(9.8) $$q_0 = -i \otimes 1 - 1 \otimes j + 1 \otimes i$$

*Тогда*

(9.9) $$\begin{aligned}s_0 &= s_1(x) - q_0 \circ p(x) = s_1(x) - (-j \otimes 1) \circ (x - i) \\ &= -jx - k + j(x - i) = -jx - k + jx + k = 0\end{aligned}$$

*Равенство*

(9.10) $$\begin{aligned}r(x) &= s_1(x) + (x \otimes 1) \circ (x - i) \\ &= (-i \otimes 1 - 1 \otimes j + 1 \otimes i + x \otimes 1) \circ (x - i)\end{aligned}$$

*является следствием равенств* (9.7), (9.9). *Мы можем упростить выражение в равенстве* (9.10)

(9.11) $$\begin{aligned}r(x) &= -i(x-i) - (x-i)j + (x-i)i + x(x-i) \\ &= (x-i)(i-j) + (x-i)(x-i) \\ &= (x-i)(x-i+i-j) \\ &= (x-i)(x-j)\end{aligned}$$

*Мы видим, что представление многочлена $r(x)$ в равенстве* (9.11) *совпадает с представление многочлена $r(x)$ в равенстве* (9.6). □

Частное от деления многочлена $r(x)$ на многочлен $p(x)$ в примере 9.2 имеет более сложную структуру чем частное в примере 9.1. Это связано с тем, что я пытаюсь записать делитель $p(x)$ справа от частного, хотя изначально многочлен $p(x)$ был левым множителем. Тем не менее, в результате несложных преобразований разложение на множители в равенстве (9.11) приняло ожидаемый вид.

$$\boxed{p(x) = x - i} \quad \boxed{r(x) = (x-i)(x-j) = x^2 - ix - xj + k}$$

В примере 9.3 я рассмотрел деление левостороннего многочлена

$$r(x) = x^2 - ix - jx - k$$

на многочлен

$$p(x) = x - i$$

Я знаю, что многочлен $r(x)$ имеет корень $x = i$, и хочу найти другой корень.

**Пример 9.3.** *Пусть*
$$p(x) = x - i$$
$$r(x) = x^2 - ix - jx - k$$

*Согласно теореме* 8.5, *положим* $q_1(x) = x \otimes 1$. *Тогда*

(9.12)
$$\begin{aligned}s_1(x) &= r(x) - (x \otimes 1) \circ (x - i) \\ &= x^2 - ix - jx - k - x(x-i) = x^2 - ix - jx - k - x^2 + xi \\ &= -ix - jx - k + xi\end{aligned}$$

*Теперь мы положим*

(9.13)
$$q_0 = -i \otimes 1 - j \otimes 1 + 1 \otimes i$$

*Тогда*

(9.14)
$$\begin{aligned}s_0 &= s_1(x) - q_0 \circ p(x) = s_1(x) - (-i \otimes 1 - j \otimes 1 + 1 \otimes i) \circ (x - i) \\ &= -ix - jx - k + xi + i(x-i) + j(x-i) - (x-i)i \\ &= -ix - jx - k + xi + ix + 1 + jx + k - xi - 1 = 0\end{aligned}$$

*Равенство*

(9.15)
$$\begin{aligned}r(x) &= s_1(x) + (x \otimes 1) \circ (x - i) = q_0 \circ (x - i) + (x \otimes 1) \circ (x - i) \\ &= (-i \otimes 1 - j \otimes 1 + 1 \otimes i + x \otimes 1) \circ (x - i)\end{aligned}$$

*является следствием равенств* (9.12), (9.14). *Мы можем упростить выражение в равенстве* (9.15)

(9.16)
$$\begin{aligned}r(x) &= -i(x-i) - j(x-i) + (x-i)i + x(x-i) \\ &= (x-i)i + (x-i-j)(x-i)\end{aligned}$$

$\square$

**Вопрос 9.4.** *Имеет ли многочлен*
$$r(x) = x^2 - ix - jx - k$$
*корень, отличный от* $x = i$? $\square$

Утверждение, что $x = i$ является корнем многочлена $r(x)$ проверяется непосредственно. Согласно примеру 9.3, мы можем представить многочлен $r(x)$ в виде

(9.17)
$$r(x) = (1 \otimes i + (x - i - j) \otimes 1) \circ (x - i)$$

Следовательно, для того чтобы многочлен $p(x)$ имел корень, отличный от $x = i$, необходимо чтобы линейное отображение

(9.18)
$$1 \otimes i + (x - i - j) \otimes 1$$

имело нетривиальное ядро и кватернион $x - i$ принадлежал ядру линейного отображения (9.18). Линейное отображение (9.18) имеет матрицу

$$(9.19) \quad \begin{pmatrix} 0 & -1 & 0 & 0 \\ 1 & 0 & 0 & 0 \\ 0 & 0 & 0 & 1 \\ 0 & 0 & -1 & 0 \end{pmatrix} + \begin{pmatrix} x^0 & -x^1+1 & -x^2+1 & -x^3 \\ x^1-1 & x^0 & -x^3 & x^2-1 \\ x^2-1 & x^3 & x^0 & -x^1+1 \\ x^3 & -x^2+1 & x^1-1 & x^0 \end{pmatrix}$$

$$= \begin{pmatrix} x^0 & -x^1 & -x^2+1 & -x^3 \\ x^1 & x^0 & -x^3 & x^2-1 \\ x^2-1 & x^3 & x^0 & -x^1+2 \\ x^3 & -x^2+1 & x^1-2 & x^0 \end{pmatrix}$$

Определить значение $x$, для которого матрица (9.19) вырождена - задача непростая. Но я надеюсь провести дополнительное исследование.

В примере 9.5 я рассмотрел деление левостороненего многочлена

$$r(x) = x^2 - ix - jx - k \qquad r(j) = -2k$$

на многочлен

$$p(x) = x - j$$

**Пример 9.5.** *Пусть*

$$p(x) = x - j$$
$$r(x) = x^2 - ix - jx - k$$

*Согласно теореме* 8.5, *положим* $q_1(x) = x \otimes 1$. *Тогда*

$$(9.20) \quad \begin{aligned} s_1(x) &= r(x) - (x \otimes 1) \circ (x - i) \\ &= x^2 - ix - jx - k - x(x - j) = x^2 - ix - jx - k - x^2 + xj \\ &= -ix - jx - k + xj \end{aligned}$$

*Теперь мы положим*

$$(9.21) \quad q_0 = -i \otimes 1 - j \otimes 1 + 1 \otimes j$$

*Тогда*

$$(9.22) \quad \begin{aligned} s_0 &= s_1(x) - q_0 \circ p(x) = s_1(x) - (-i \otimes 1 - j \otimes 1 + 1 \otimes j) \circ (x - j) \\ &= -ix - jx - k + xj + i(x - j) + j(x - j) - (x - j)j \\ &= -ix - jx - k + xj + ix - k + jx + 1 - xj - 1 = -2k \end{aligned}$$

*Равенство*

$$(9.23) \quad \begin{aligned} r(x) &= s_1(x) + (x \otimes 1) \circ (x - j) \\ &= s_0 + q_0 \circ (x - j) + (x \otimes 1) \circ (x - j) \\ &= -2k + (-i \otimes 1 - j \otimes 1 + 1 \otimes j + x \otimes 1) \circ (x - j) \end{aligned}$$

*является следствием равенств* (9.20), (9.22). *Мы можем упростить выражение в равенстве* (9.23)

$$(9.24) \quad r(x) = -2k + (1 \otimes j + (x - i - j) \otimes 1) \circ (x - j)$$

$\square$

## 10. Новые вопросы

В примере 10.1 я рассмотрел многочлен $r(x)$ степени 3 для того, чтобы найти частное от деления многочлен $r(x)$ на два многочлена. Сперва я найду частное от деления многочлена

$$r(x) = (x - j)(x - k)(x - j - k)$$

на многочлен

$$p(x) = x - k$$

**Пример 10.1.** *Пусть*

$$p(x) = x - k$$

$$\begin{aligned}
r(x) &= (x - j)(x - k)(x - j - k) = (x^2 - jx - xk + i)(x - j - k) \\
&= x^2(x - j - k) - jx(x - j - k) - xk(x - j - k) + i(x - j - k) \\
&= x^3 - x^2 j - x^2 k - jx^2 + jxj + jxk - xkx + xkj - x + ix - k + j \\
&= x^3 \\
&\quad + (-1 \otimes 1 \otimes j - 1 \otimes 1 \otimes k - j \otimes 1 \otimes 1 - 1 \otimes k \otimes 1) \circ x^2 \\
&\quad + (j \otimes j + j \otimes k - 1 \otimes i - 1 \otimes 1 + i \otimes 1) \circ x \\
&\quad - k + j
\end{aligned}$$

*Согласно теореме* 8.5*, положим* $q_2(x) = (x^2) \otimes 1$. *Тогда*

$$\begin{aligned}
s_2(x) &= r(x) - q_2(x) \circ (x - k) = r(x) - (x^2 \otimes 1) \circ (x - k) \\
&= x^3 - x^2 j - x^2 k - jx^2 + jxj + jxk - xkx + xkj - x + ix - k + j \\
&\quad - x^3 + x^2 k \\
&= -x^2 j - jx^2 + jxj + jxk - xkx + xkj - x + ix - k + j
\end{aligned} \tag{10.1}$$

*Согласно теореме* 8.5*, положим*

$$q_1(x) = -x \otimes j - jx \otimes 1 - xk \otimes 1 \tag{10.2}$$

*Равенство*

$$\begin{aligned}
s_1(x) &= s_2(x) - q_1(x) \circ (x - k) \\
&= -x^2 j - jx^2 + jxj + jxk - xkx + xkj - x + ix - k + j \\
&\quad - (-x \otimes j - jx \otimes 1 - xk \otimes 1) \circ (x - k) \\
&= -x^2 j - jx^2 + jxj + jxk - xkx + xkj - x + ix - k + j \\
&\quad + x(x - k)j + jx(x - k) + xk(x - k) \\
&= -x^2 j - jx^2 + jxj + jxk - xkx + xkj - x + ix - k + j \\
&\quad + x^2 j - xkj + jx^2 - jxk + xkx - xkk \\
&= jxj + ix - k + j
\end{aligned} \tag{10.3}$$

*является следствием равенств* (10.1), (10.2). *Согласно теореме* 8.5*, положим*

$$q_0(x) = j \otimes j + i \otimes 1 \tag{10.4}$$

*Равенство*

$$\begin{aligned}s_0 &= s_1(x) - q_0(x) \circ (x-k)\\ &= jxj + ix - k + j - (j \otimes j + i \otimes 1) \circ (x-k)\\ &= jxj + ix - k + j - j(x-k)j - i(x-k)\\ &= jxj + ix - k + j - jxj + jkj - ix + ik\\ &= -k + j + jkj + ik = 0\end{aligned}$$
(10.5)

*является следствием равенств* (10.3), (10.4). *Равенство*

$$\begin{aligned}r(x) &= s_0 + (q_0 + q_1(x) + q_2(x)) \circ (x-k)\\ &= (j \otimes j + i \otimes 1 - x \otimes j - jx \otimes 1 - xk \otimes 1 + x^2 \otimes 1) \circ (x-k)\end{aligned}$$
(10.6)

*является следствием равенств* (10.1), (10.2), (10.4), (10.5).　□

Частное

(10.7) $\quad q(x) = j \otimes j + i \otimes 1 - x \otimes j - jx \otimes 1 - xk \otimes 1 + x^2 \otimes 1$

от деления многочлена

$$r(x) = (x-j)(x-k)(x-j-k)$$

на многочлен

$$p(x) = x - k$$

является тензором зависящим от $x$, а не многочленом. Поэтому мы не можем применить теорему 8.5.

**Вопрос 10.2.** *Нетрудно убедиться, что*

$$\begin{aligned}q(j) &= (j \otimes j + i \otimes 1 - x \otimes j - jx \otimes 1 - xk \otimes 1 + x^2 \otimes 1)(j)\\ &= j \otimes j + i \otimes 1 - j \otimes j - j^2 \otimes 1 - jk \otimes 1 + j^2 \otimes 1 = 0\end{aligned}$$
(10.8)

*Поэтому возникает вопрос как выделить в тензоре*

(10.9) $\quad j \otimes j + i \otimes 1 - x \otimes j - jx \otimes 1 - xk \otimes 1 + x^2 \otimes 1$

*множитель* $x - j$.　□

**Пример 10.3.** *Один из возможных ответов на вопрос* 10.2 *состоит в преобразовании*

$$\begin{aligned}&j \otimes j + i \otimes 1 - x \otimes j - jx \otimes 1 - xk \otimes 1 + x^2 \otimes 1\\ =&j \otimes j - x \otimes j + i \otimes 1 - jx \otimes 1 - xk \otimes 1 + x^2 \otimes 1\\ =&-(x-j) \otimes j + (i - jx - xk + x^2) \otimes 1\end{aligned}$$
(10.10)

*Рассмотрим многочлен*

(10.11) $\quad\quad\quad\quad\quad\quad i - jx - xk + x^2$

*Согласно теореме* 8.5, *положим* $q_1(x) = x \otimes 1$. *Тогда*

$$\begin{aligned}s_1(x) &= r(x) - q_1(x) \circ (x-j) = r(x) - (x \otimes 1) \circ (x-j)\\ &= i - jx - xk + x^2 - x^2 + xj\\ &= i - jx - xk + xj\end{aligned}$$
(10.12)

*Теперь мы положим*

(10.13) $$q_0 = -j \otimes 1 - 1 \otimes k + 1 \otimes j$$

*Тогда*

(10.14) $$\begin{aligned} s_0 &= s_1(x) - q_0 \circ p(x) = s_1(x) - (-j \otimes 1 - 1 \otimes k + 1 \otimes j) \circ (x - j) \\ &= i - jx - xk + xj + j(x-j) + (x-j)k - (x-j)j \\ &= i - jx - xk + xj + jx - j^2 + xk - jk - xj + j^2 \\ &= i - jk = 0 \end{aligned}$$

*Равенство*

(10.15) $$\begin{aligned} r(x) &= s_1(x) + (x \otimes 1) \circ (x - j) \\ &= (-j \otimes 1 - 1 \otimes k + 1 \otimes j + x \otimes 1) \circ (x - j) \end{aligned}$$

*является следствием равенств* (10.12), (10.14). *Равенство*

(10.16) $$\begin{aligned} &j \otimes j + i \otimes 1 - x \otimes j - jx \otimes 1 - xk \otimes 1 + x^2 \otimes 1 \\ &= -(x-j) \otimes j + ((-j \otimes 1 - 1 \otimes k + 1 \otimes j + x \otimes 1) \circ (x - j)) \otimes 1 \\ &= ((-1 \otimes j - j \otimes 1 - 1 \otimes k + 1 \otimes j + x \otimes 1) \circ (x - j)) \otimes 1 \\ &= ((-j \otimes 1 - 1 \otimes k + x \otimes 1) \circ (x - j)) \otimes 1 \end{aligned}$$

*является следствием равенств* (10.10), (10.15). □

Мы видим, что частное в примере 10.3 является $H^{3\otimes}$-числом.

**Пример 10.4.** *Равенство*

(10.17) $$\begin{aligned} r(x) &= (((-j \otimes 1 - 1 \otimes k + 1 \otimes j + x \otimes 1) \circ (x - j)) \otimes 1) \circ (x - k) \\ &= (-j \otimes 1 \otimes 1 - 1 \otimes k \otimes 1 + 1 \otimes j \otimes 1 + x \otimes 1 \otimes 1) \circ (x - j, x - k) \end{aligned}$$

*является следствием равенств* (10.6), (10.16). □

**Вопрос 10.5.** *Равенство* (10.17) *является ответом на вопрос* 10.2. *Мы видим, что разложение многочлена* $r(x)$ *является билинейным отображением многочленов* $x - j$, $x - k$. *Однако этот ответ не полный. Согласно замечанию* 2.6, *представление* (8.5) *многочлена не единственно. Остаётся открытым вопрос о выборе множителей и их порядке. В то же время существует связь между различными представлениями многочлена в виде произведения множителей.* □

Если многочлен $p(x)$ имеет вид $p(x) = p_1 \circ x + p_0$ где $p_1$ - невырожденный тензор, то деление на многочлен $p(x)$ можно свести к теореме 8.5, если мы воспользуемся равенством

(10.18) $$p(x) = p_1 \circ (x + p_1^{-1} \circ p_0)$$

**Вопрос 10.6.** *Многочлен*

(10.19) $$p_1(x) = ix - xi - 1$$

*является делителем многочлена*

$$\begin{aligned}
p(x) &= (x-j)(ix - xi - 1)(x-k) \\
&= ((x-j)ix - (x-j)xi - x + j)(x-k) \\
&= xix(x-k) - jix(x-k) - x^2 i(x-k) \\
&\quad + jxi(x-k) - x(x-k) + j(x-k) \\
&= xix^2 - xixk + kx^2 - kxk - x^2 ix + x^2 ik \\
&\quad + jxix - jxik - x^2 + xk + jx - jk \\
&= (1 \otimes i \otimes 1 \otimes 1 - 1 \otimes 1 \otimes i \otimes 1) \circ x^3 \\
&\quad + (j \otimes i \otimes 1 - 1 \otimes i \otimes k + k \otimes 1 \otimes 1 - 1 \otimes 1 \otimes j - 1 \otimes 1 \otimes 1) \circ x^2 \\
&\quad + (1 \otimes k + j \otimes (1+j)) \circ x + 1 - i
\end{aligned} \tag{10.20}$$

*Однако мы не можем воспользоваться теоремой 8.5 для нахождения частного, так как тензор*

$$i \otimes 1 - 1 \otimes i \in H^{2\otimes}$$

*вырожден. Поэтому необходимо дополнительное исследование, чтобы решить эту задачу.* □

## 11. Деление многочленов в неассоциативной алгебре

**Замечание 11.1.** *В неассоциативной алгебре скобки определяют в каком порядке мы выполняем умножение. Поэтому одночлен может иметь следующий вид*

$$(p(x)q(x))r(x) \tag{11.1}$$

$$p(x)(q(x)r(x)) \tag{11.2}$$

*где $p$, $q$, $r$ также одночлены. Это накладывает ограничение на возможность деления многочленов.* □

В примере 11.2 я рассмотрел многочлен $r(x)$ степени 3 в алгебре октонионов. Так как порядок множителей в произведении важен для определения структуры тензора, я буду явно указывать скобки даже в тех случаях, когда порядок очевиден. Например, поскольку, вообще говоря,

$$((a\otimes)b) \circ x = (ax)b \neq a(xb) = (a(\otimes b)) \circ x \tag{11.3}$$

то выражения $(a\otimes)b$ и $a(\otimes b)$ описывают различные линейные отображения.

**Пример 11.2.** *Пусть*

$$\begin{aligned}
p(x) &= x - k \\
r(x) &= ((x-j)(x-k))(x-jl) = (x^2 - jx - xk + i)(x-jl) \\
&= (x^2)(x-jl) - (jx)(x-jl) - (xk)(x-jl) + i(x-jl) \\
&= (x^2)x - (x^2)jl - (jx)x - (xk)x + (jx)jl + (xk)jl + ix - i\,jl
\end{aligned} \tag{11.4}$$

*Согласно теореме* 8.5, *положим* $q_2(x) = (x^2)(\otimes 1)$. *Тогда*

$$\begin{aligned}
s_2(x) &= r(x) - q_2(x) \circ (x-k) = r(x) - ((x^2)(\otimes 1)) \circ (x-k) \\
&= (x^2)x - (x^2)jl - (jx)x - (xk)x + (jx)jl + (xk)jl + ix + kl \\
&\quad - (x^2)x + (x^2)k \\
&= (x^2)k - (x^2)jl - (jx)x - (xk)x + (jx)jl + (xk)jl + ix + kl
\end{aligned}$$
(11.5)

*Согласно теореме* 8.5, *положим*

(11.6) $$q_1(x) = (x\otimes)k - (x\otimes)jl - (jx)(\otimes 1) - (xk)(\otimes 1)$$

*Равенство*

$$\begin{aligned}
s_1(x) &= s_2(x) - q_1(x) \circ (x-k) \\
&= (x^2)k - (x^2)jl - (jx)x - (xk)x + (jx)jl + (xk)jl + ix + kl \\
&\quad - ((x\otimes)k - (x\otimes)jl - (jx)(\otimes 1) - (xk)(\otimes 1)) \circ (x-k) \\
&= (x^2)k - (x^2)jl - (jx)x - (xk)x + (jx)jl + (xk)jl + ix + kl \\
&\quad - (x(x-k))k + (x(x-k))jl + (jx)(x-k) + (xk)(x-k) \\
&= (x^2)k - (x^2)jl - (jx)x - (xk)x + (jx)jl + (xk)jl + ix + kl \\
&\quad - (x^2)k + (xk)k + (x^2)jl - (xk)jl + (jx)x - (jx)k + (xk)x - (xk)k \\
&= (jx)(jl-k) + ix + kl
\end{aligned}$$
(11.7)

*является следствием равенств* (11.5), (11.6). *Согласно теореме* 8.5, *положим*

(11.8) $$q_0(x) = (j\otimes)(jl-k) + i(\otimes 1)$$

*Равенство*

$$\begin{aligned}
s_0 &= s_1(x) - q_0(x) \circ (x-k) \\
&= (jx)(jl-k) + ix + kl - ((j \otimes (jl-k)) + i(\otimes 1)) \circ (x-k) \\
&= (jx)(jl-k) + ix + kl - (j(x-k))(jl-k) - i(x-k) \\
&= (jx)(jl-k) + ix + kl - (jx - jk)(jl-k) - ix + ik \\
&= (jx)(jl-k) + ix + kl - (jx)(jl-k) + i(jl-k) - ix - j \\
&= kl + ijl - ik - j = 0
\end{aligned}$$
(11.9)

*является следствием равенств* (11.7), (11.8). *Равенство*

$$\begin{aligned}
r(x) &= s_0 + (q_0 + q_1(x) + q_2(x)) \circ (x-k) \\
&= ((j\otimes)(jl-k) + i(\otimes 1) \\
&\quad + (x\otimes)k - (x\otimes)jl - (jx)(\otimes 1) - (xk)(\otimes 1) + (x^2)(\otimes 1)) \circ (x-k)
\end{aligned}$$
(11.10)

*является следствием равенств* (11.5), (11.6), (11.8), (11.9). □

**Пример 11.3.** *Частное от деления многочлена* (11.4) *на многочлен*
$$p(x) = x - k$$
*имеет вид*

$$\begin{aligned}
q(x) &= (j\otimes)(jl-k) + i(\otimes 1) \\
&\quad + (x\otimes)k - (x\otimes)jl - (jx)(\otimes 1) - (xk)(\otimes 1) + (x^2)(\otimes 1) \\
&= ((x-j)\otimes)(k-jl) + (x^2 - jx - xk + i)(\otimes 1)
\end{aligned}$$
(11.11)

(11.4) $\boxed{r(x) = ((x-j)(x-k))(x-jl)}$

Первое слагаемое тензора $q(x)$ имеет множитель $x-j$. Чтобы проверить, имеет ли второе слагаемое множитель $x-j$, рассмотрим многочлен

(11.12) $$r_1(x) = x^2 - jx - xk + i$$

Согласно теореме 8.5, положим $q_1(x) = x(\otimes 1)$. Тогда

(11.13) $$\begin{aligned} s_1(x) &= r_1(x) - (x(\otimes 1)) \circ (x-j) \\ &= x^2 - jx - xk + i - x(x-j) = x^2 - jx - xk + i - x^2 + xj \\ &= -jx - xk + xj + i \end{aligned}$$

Теперь мы положим

(11.14) $$q_0 = -j(\otimes 1) - (1\otimes)k + (1\otimes)j$$

Тогда

(11.15) $$\begin{aligned} s_0 &= s_1(x) - q_0 \circ p(x) = s_1(x) - (-j(\otimes 1) - (1\otimes)k + (1\otimes)j) \circ (x-j) \\ &= -jx - xk + xj + i + j(x-j) + (x-j)k - (x-j)j \\ &= -jx - xk + xj + i + jx + 1 + xk - i - xj - 1 = 0 \end{aligned}$$

Равенство

(11.16) $$\begin{aligned} r_1(x) &= s_1(x) + (x(\otimes 1)) \circ (x-j) \\ &= (-j(\otimes 1) - (1\otimes)k + (1\otimes)j + x(\otimes 1)) \circ (x-j) \\ &= ((1\otimes)(j-k) + (x-j)(\otimes 1)) \circ (x-j) \end{aligned}$$

является следствием равенств (11.13), (11.15). $\square$

Из равенств (11.11), (11.16) следует, что частное от деления многочлена (11.4) на многочлен

$$p(x) = x - k$$

имеет вид

(11.17) $$\begin{aligned} q(x) &= = ((x-j)\otimes)(k-jl) \\ &+ (((1\otimes)(j-k) + (x-j)(\otimes 1)) \circ (x-j))(\otimes 1) \\ &= (((1\otimes_1)1\otimes)(k-jl) \\ &+ ((1\otimes_1)(j-k) + (x-j)(\otimes_1 1))(\otimes 1)) \circ (x-j) \\ &= (((1\otimes_1)1\otimes)(k-jl) \\ &+ ((1\otimes_1)(j-k))(\otimes 1) + ((x-j)(\otimes_1 1))(\otimes 1)) \\ &\circ (x-j) \end{aligned}$$

(11.4) $\boxed{r(x) = ((x-j)(x-k))(x-jl)}$

Следовательно, мы можем представить многочлен (11.4) в виде

(11.18) $$\begin{aligned} r(x) &= (((1\otimes_1)1\otimes_2)(k-jl) \\ &+ ((1\otimes_1)(j-k))(\otimes_2 1) + ((x-j)(\otimes_1 1))(\otimes_2 1)) \\ &\circ (x-j, x-k) \end{aligned}$$

## 12. Список литературы

## 13. Предметный указатель

## 14. Специальные символы и обозначения